\documentclass[a4wide,12pt]{article}
\usepackage{amsmath,multicol,graphicx,pdfpages}

\begin{document}

\title{An investigation of the comparative efficiency of the different methods in which $\pi$ is calculated}
\author{Nouri Al-Othman}
\maketitle
\newpage

\abstract
An investigation of the comparative efficiency of the different methods in which $\pi$ is calculated”. This thesis will compare and contrast five different methods in calculating $\pi$ by first deriving the various proofs to each method and then creating graphical displays and tables with percentage errors on each to allow a thorough comparison between each method. It is important to realize that usually people refer to $\pi$ as 3.142 and use this shortened value in everyday calculations. However, most people in society do not know that in some of NASA’s calculations, π must be calculated to ”15 decimal places (3.141592653589793) for calculations involving the Space Station Guidance Navigation and Control (GNC) code” [9]. Furthermore, some people use the world renown constant everyday yet have no idea where it comes from.
\\ \\
This thesis allows a thorough exploration of the different methods of calculating pi by using Mathematica to provide the necessary iterations needed and Microsoft excel to plot the various graphs for convergence comparisons. Finally, the conclusions that have come from this thesis demonstrate that methods in which oscillations between the positive and negative region of the line y = $\pi$ slow down the rate of convergence and out of the five methods investigated, The infinite series formulated from $\zeta(8)$ proves to be one of the fastest infinite series today to converge towards $\pi$.

\newpage

\tableofcontents
\newpage

\section{introduction}
In today's society, $\pi$ is used everyday by almost everyone in some way or another. Wether it be a mathematician doing complex integrals, an engineer calculating the volume of spherical shapes or a clock maker designing a new pendulum for a new clock. It is clear that $\pi$ underlies most of the modern day achievements and developments. However, the reality is that only a small percentage of the population who use this constant have any idea where it comes from.  \\ \\

In the beginning, Pi was first investigated using a trigonometric method in "250BC by the famous mathematician Archimedes" \cite{cite_key7} who investigated interior and exterior polygons around a circle creating both upper and lower bounds for a value of $\pi$ which were then averaged to give a value correctly to 7 decimal places \cite{cite_key7}.  Moving swiftly through the decades, another mathematician named Vieta in 1579 used nested radicals to further find a product of nested radicals which would converge to $\pi$ and gave a value of $\pi$ correct to 9 deicmal places. Lastly, with the development of calculus, rapidly converging series were developed. A clear example of this is when newton in 1748 developed his arcsine expansion with the Taylor Series which gave a value of $\pi$ which "only needed 22 terms for 16 decimal places for $\pi$" \cite{cite_key5}.  \\ \\

However, with the mathematical knowledge we have today, there are endless ways to arrive at this world renown constant of $\pi$. From the simplest application of this constant in the formula $\pi r^2$ to attain a value for the area of a circle to the more complex application towards the other end of the spectrum in calculus. My research question: ``An investigation of the comparative efficiency of the different methods in which $\pi$ is calculated" aims to distinguish and set aside the different methods of calculating and approaching the value of this constant by contrasting how efficient a variety of different numerical methods are in calculating this constant correct to 15 decimal places. This will be achieved by exploring different methods in which $\pi$ is calculated and contrasting their rates of convergence by increasing the amount of terms $n$ and further observing the rate at which each method converges towards $\pi$. A series of graphical displays and Data tables will be used to calculate and compare the convergence towards $\pi$ and where there is an error by comparing it to a "True" value.
\newpage

\section{The Riemann-Zeta function, Parseval's identity and the Fourier Series to generate $\zeta (x)$}
\subsection{Preface list of data and proof's used in the calculations}
\subsubsection{The Fourier Series}

\begin{equation}
 (s_N f)(x) = \frac{a_0}{2} +\sum_{n=1}^N[a_n cos(nx) + b_n sin(nx)]  \mbox {   for   } N \in Z^+ \cite{cite_key10} \nonumber
\end{equation}

\subsubsection{Parseval's Identity}
\begin{equation}
\sum_{\infty}^{-\infty} |{c_n}|= \frac{1}{\pi} \int_{-\pi}^{\pi} |f(x)|^2  \cite{cite_key6} \nonumber \mbox {     where} \frac{1}{\pi} \int^{\pi}_{-\pi} f^2 (x) dx = \frac{{a_0}^2}{2} + \sum^{\infty}_{k=1} | {a_k}^2 + {b_k}^2 |
\end{equation}

\subsubsection{Proof of orthongonality in the Fourier Series}
First consider
\begin{equation}
f(x) = \frac{a_0}{2} + \sum_{n=1}^N a_n cos(nx) + \sum_{n=1}^N b_n sin(nx) \nonumber
\end{equation}
Case (I) Note: $sin[(n+m)x] + sin[(n-m)x] = 2sin(nx)cos(mx)$

\begin{eqnarray}
\int_{-\pi}^{\pi} sin(nx)cos(mx)dx &=& \frac{1}{2} \int_{-\pi}^{\pi} 2sin(nx)cos(mx)dx \nonumber \\
 &=& \frac{1}{2} \int_{-\pi}^{\pi} sin[(n+m)x]+sin[(n-m)x]dx \nonumber \\
 &=& \frac{1}{2} \left[\frac{-1}{(n+m)} cos[(n+m)x] + \frac{1}{(n-m)} cos[(n-m)x]\right]_{-\pi}^{\pi} \nonumber \\
 &=& \left[\frac{-1}{2(n+m)} cos[(n+m)x] + \frac{1}{2(n-m)} cos[(n-m)x]\right]_{-\pi}^{\pi} \nonumber\\
 &=& 0 \nonumber 
\end{eqnarray}

Case (II) Note: $cos[(n+m)x] + cos[(n-m)x] = 2cos(nx)cos(mx)$

\begin{eqnarray}
\int_{-\pi}^{\pi} cos(nx)cos(mx)dx &=& \frac{1}{2} \int_{-\pi}^{\pi} 2cos(nx)cos(mx)dx \nonumber \\
 &=& \frac{1}{2} \int_{-\pi}^{\pi} cos[(n+m)x]+cos[(n-m)x]dx \nonumber \\
 &=& \frac{1}{2} \left[\frac{1}{(n+m)} sin[(n+m)x] + \frac{1}{(n-m)} sin[(n-m)x]\right]_{-\pi}^{\pi} \mbox{           Note: Only if $n \neq m$} \nonumber\\
 &=& \left[\frac{1}{2(n+m)} sin[(n+m)x] + \frac{1}{2(n-m)} sin[(n-m)x]\right]_{-\pi}^{\pi} \nonumber
\end{eqnarray}

Firstly, consider the case where $n \neq m$ \\
\begin{eqnarray}
\int_{-\pi}^{\pi} cos(nx)cos(mx) = 0 \nonumber
\end{eqnarray}

Secondly, consider the case where $n = m$ \\
\begin{eqnarray}
\int^{\pi}_{-\pi} cos(nx)cos(nx) dx &=& \int^{\pi}_{-\pi} cos^2 (nx) dx = \frac{1}{2} \int^{\pi}_{-\pi} (cos(2nx) +1) dx = \frac{1}{2} \left[\frac{1}{2n} sin(2nx) + x \right]^{\pi}_{-\pi} \nonumber \\
&=& \frac{1}{2} \left[0 + 2\pi \right] \nonumber \\
&=& \pi \nonumber
\end{eqnarray}

Case (III) note: $-cos[(n+m)x] + cos[(n-m)x]$ = $2sin(nx)sin(mx)$ \\
\begin{eqnarray} 
\int^{\pi}_{-\pi} sin(mx)sin(nx) &=& \frac{1}{2} \int^\pi_-\pi 2sin(mx)sin(nx) dx \nonumber \\ 
&=& \frac{1}{2} \int^{\pi}_{-\pi} cos[(n-m)x] - cos[(n-m)x] dx \nonumber \\
&=& \frac{1}{2} \left[\frac{1}{(n-m)} sin[(n-m)x] - \frac{1}{(n+m)} sin[(n+m)x]\right]_{-\pi}^{\pi} \mbox{ for $m \neq n$} \nonumber
\end{eqnarray}

Case (i): If $n \neq m$  then, 
\begin{eqnarray}
\int_{-\pi}^{\pi} sin(mx)sin(nx) dx = 0 \nonumber 
\end{eqnarray}

Case (ii): If $n = m$ \ then,
\begin{eqnarray}
\int^{\pi}_{-\pi} sin(mx)sin(nx) &=& \int^{\pi}_{-\pi} sin^2 (nx) dx \nonumber \\
&=& \frac{1}{2} \int^{\pi}_{-\pi} (1-cos(2nx)) dx \nonumber \\
&=&\frac{1}{2} \left[x - \frac{1}{2n} sin(2nx) \right]^{\pi}_{-\pi} \nonumber \\ 
&=& \frac{1}{2} \left[2\pi - 0\right] \nonumber \\
&=& \pi \nonumber
\end{eqnarray}

Therefore, we can deduce that in general:

\begin{eqnarray}
\int_{-\pi}^{\pi} sin(nx)cos(mx) dx &=& 0  \nonumber \\ \nonumber \\
\int_{-\pi}^{\pi} sin(nx)sin(mx) dx &=& \nonumber
\left\{ 
\begin{array}{rl}
0 &\mbox{ if $m \neq n$} \\ 
\pi &\mbox{ if $m = n$} 
\end{array}
 \right.
\rightarrow \frac{1}{\pi}  \int_{-\pi}^{\pi} sin(nx)sin(mx) dx =
\left\{ 
\begin{array}{rl}
0 &\mbox{ if $m \neq n$} \\
1 &\mbox{ if $m = n$} \nonumber
\end{array} \right.
\end{eqnarray}

\begin{eqnarray}
\int_{-\pi}^{\pi} cos(nx)cos(mx) dx &=& \nonumber
\left\{ 
\begin{array}{rl}
0 &\mbox{ if $m \neq n$} \\ 
\pi &\mbox{ if $m = n$} 
\end{array}
 \right.
\rightarrow \frac{1}{\pi}  \int_{-\pi}^{\pi} cos(nx)cos(mx) dx =
\left\{ 
\begin{array}{rl}
0 &\mbox{ if $m \neq n$} \\
1 &\mbox{ if $m = n$} \nonumber
\end{array} \right.
\end{eqnarray}

\newpage
\subsubsection{The extraction of $a_n , b_n , a_0, b_0$ from the Fourier series}

From the Fourier series, by applying integration on f(x) directly and realizing that the terms consisting of $sin(nx)$ and $cos(nx)$ are equal to zero, we are able to attain a value for the constant term $a_0$.

\begin{eqnarray}
\int^{\pi}_{-\pi} f(x) dx &=& \int^{\pi}_{-\pi} 
\left(\frac{a_0}{2}\right) dx + \sum_{n=1}^{\infty} a_n \int^{\pi}_{-\pi} cos(nx) dx +  \sum_{n=1}^{\infty} b_n \int^{\pi}_{-\pi} sin(nx) dx \nonumber \\
&=& \left(\frac{a_0}{2}\right) \int_{-\pi}^{\pi} \frac{dx}{2} \nonumber \\
&=& a_0 \left[\frac{x}{2} \right]^{\pi}_{-\pi} \nonumber \\
&=& a_0 \pi \nonumber \\ \nonumber \\
a_0 \pi &=& \int^{\pi}_{-\pi} f(x) dx \rightarrow a_0 = \frac{1}{\pi} \int^{\pi}_{-\pi} f(x) dx \nonumber
\end{eqnarray}

To find a value of $a_n$, multiply all the terms in the Fourier series by cos(nx) as $a_n$ will only exist where f(x) is even because the term in $a_n$ will be the only one to be non-zero following the integration over the range from $\pi$ to $-\pi$ \nonumber

\begin{eqnarray}
\int^{\pi}_{-\pi} f(x)cos(mx) dx &=& \int^{\pi}_{-\pi} 
\left(\frac{a_0}{2}\right) cos(mx) dx  + \sum_{m=1}^{\infty} a_m \int^{\pi}_{-\pi} cos(nx)cos(mx) dx +  \sum_{m=1}^{\infty} b_m \int^{\pi}_{-\pi} sin(nx)cos(mx) dx \nonumber \\
&=& \sum_{m=1}^{\infty} a_m \int^{\pi}_{-\pi} cos(nx)cos(mx) dx \nonumber \\
&=& a_m \int^{\pi}_{-\pi} cos(nx)cos(mx)  dx \nonumber
\end{eqnarray}

Note that $\delta_{mn} = 0$ where $m \neq n$ and $\delta_{mn} = 1 $ where $m = n$ \\ \\

Now note that there only exists $a_n$ where $m = n$

\begin{eqnarray}
\int_{-\pi}^{\pi} f(x)cos(nx) dx = a_n\pi \rightarrow a_n =\frac {1}{\pi} \int^{\pi}_{-\pi} f(x)cos(nx) dx \nonumber
\end{eqnarray}

In addition, when we multiply all terms by sin(nx) we can yield a resut respecitvly for $b_n$

\begin{eqnarray}
b_n = \frac{1}{\pi} \int^{\pi}_{-\pi} f(x)sin(nx)dx \nonumber
\end{eqnarray}

Also, note that in all cases dealt with $b_0 = 0$ and $b_n$ is only defined for $n \geq 1$

\newpage
\subsection{Generating a value for $\zeta$ (2)}

Consider the simplest odd function where $f(x) = x$ \\
Therefore, since $f(x)$ is an odd function then $a_n = 0$ where $\forall n \in Z^+$\nonumber \\

Now consider $b_n$,

\begin{eqnarray}
b_n = \frac{1}{\pi} \int^{\pi}_{-\pi} xsin(nx) dx &=& \frac{1}{\pi} \left(\left[\frac{-xcos(nx)}{n}\right]^{\pi}_{-\pi} -\int^{\pi}_{-\pi} \frac{-cos(nx)}{n} dx\right) \nonumber \\
 &=& \frac{1}{\pi} \left[\left[\frac{sin(nx)}{n^2}\right]^{\pi}_{-\pi} - \left[\frac{xcos(nx)}{n}\right]^{\pi}_{-\pi} \right] \nonumber \\
&=& \frac{1}{\pi} \left[\frac{-xcos(nx)}{n}\right]^{\pi}_{-\pi} \nonumber \\
&=& \frac{1}{\pi} \left[\frac{-\pi cos(n\pi)}{n} - \frac{\pi cos(-n\pi)}{n} \right] \nonumber \\
&=& \frac{1}{\pi} \left[\frac{2\pi cos(n{\pi})}{n}\right] \nonumber \\
&=& \frac{-2cos(n{\pi})}{n} \nonumber \\
&=& \frac{-2(-1)^n}{n} \nonumber
\end{eqnarray}

Now, by substituting $b_n$ back into the Fourier Series we obtain:

\begin{eqnarray}
x = \sum_{n=1}^{\infty} \frac{-2(-1)^n}{n} sin(nx)     \mbox{     for     }    0\leq x\leq \pi  \nonumber  
\end{eqnarray}

Then, consider Parsevals Identity where: \\
\begin{eqnarray}
c_n = a_n + ib_n \rightarrow |c_n|^2 &=& a_n^2 + b_n^2  \mbox{    (Note: $a_n = 0$ where $f(x)$ is odd)    } \nonumber \\
\sum_{n=-\infty}^{\infty} |c_n|^2 &=& \frac{1}{\pi} \int^{\pi}_{-\pi} |f(x)|^2  dx \nonumber \\
\sum_{n=-\infty}^{\infty} |b_n|^2 &=& \frac{1}{\pi} \int^{\pi}_{-\pi} x^2  dx \nonumber
\end{eqnarray}

Firstly, consider the left hand side of Parseval's Identity:
\begin{eqnarray}
\sum_{n=1}^{\infty} |b_n|^2 &=& \sum^{\infty}_{n=1} \left[\frac{-2(-1)^n}{n}\right]^2 \nonumber \\
&=& \sum_{n=1}^{\infty} \frac{4(-1)^{2n}}{n^2} \nonumber \\
&=& 4 \sum_{n=1}^{\infty} \frac{1}{n^2} \nonumber \\
&=& 4 \zeta (2) \nonumber
\end{eqnarray}

Secondly, consider the right hand side of Parseval's Identity:
\begin{eqnarray}
\frac{1}{\pi} \int^{\pi}_{-\pi} x^2 &=& \frac{1}{\pi} \left[\frac{2\pi^3}{3}\right] \nonumber \\
&=& \frac{2\pi^2}{3} \nonumber
\end{eqnarray}

Lastly, equate the left hand side to the right hand side of Parseval's Identity
\begin{eqnarray}
4\zeta (2) &=& \frac{2\pi^2}{3} \nonumber \\
\zeta (2) &=& \frac{2\pi^2}{12} \nonumber \\
\zeta (2) &=& \frac{\pi^2}{6} \nonumber
\end{eqnarray}

Now, using the definition of the zeta function we then see that

\begin{eqnarray}
\zeta (2) = \sum^{\infty}_{k=1} \frac{1}{k^2} = \frac{1}{1^2} + \frac{1}{2^2} +\frac{1}{3^2} + ... \nonumber
\end{eqnarray}

Multiplying by 6 and square rooting we attain:

\begin{eqnarray}
{\pi^2} &=& 6\left[\frac{1}{1^2} + \frac{1}{2^2} +\frac{1}{3^2} + ... \right] \nonumber \\ \nonumber \\
\pi &=& \sqrt{6\left[\frac{1}{1^2} + \frac{1}{2^2} +\frac{1}{3^2} + ... \right] \nonumber } \nonumber
\end{eqnarray}

\newpage
\subsection{Generating a value for $\zeta$ (4)}

Consider the simplist even function where $f(x) = x^2$ \\
Therefore, since f(x) is an even function then $b_n = 0$ where $\forall n \in Z^+$ \nonumber \\

First consider $a_0$,
\begin{eqnarray}
a_0 = \frac{1}{\pi} \int^{\pi}_{-\pi} x^2 dx = \frac{2}{\pi} \int^{\pi}_{0} x^2 dx = \frac{2\pi^2}{3} \nonumber
\end{eqnarray}

Then consider $a_n$,
\begin{eqnarray}
a_n = \frac{1}{\pi} \int^{\pi}_{-\pi} x^2 cos(nx) dx &=& \frac{2}{\pi} \int^{\pi}_{0} x^2 cos(nx) dx\nonumber \\
&=& \frac{2}{\pi} \left[\left[\frac{x^2 sin(nx)}{n}\right]^{\pi}_{0} - \int^{\pi}_{0} \frac{2xsin(nx)}{n} dx\right] \nonumber \\
&=& \frac{2}{\pi} \left[\left[\frac{2x cos(nx)}{n^2}\right]^{\pi}_{0} - \int^{\pi}_{0} \frac{2cos(nx)}{n^2} dx\right] \nonumber \\
&=& \frac{4\pi cos(nx)}{n^2 \pi} \nonumber \\
&=& \frac{4(-1)^n}{n^2 } \nonumber
\end{eqnarray}

Now, by substituting $a_n$ and $a_0$ back into the Fourier Series we obtain:

\begin{eqnarray}
x = \frac{\pi^2}{3} + \sum_{n=1}^{\infty} \frac{4(-1)^n}{n^2} cos(nx)     \mbox{     for     }    -\pi\leq x\leq \pi  \nonumber  
\end{eqnarray}

Then, consider Parsevals Identity where: \\
\begin{eqnarray}
c_n = a_n + ib_n \rightarrow |c_n|^2 &=& a_n^2 + b_n^2  \mbox{    (Note: $b_n = 0$ where $f(x)$ is even)    } \nonumber \\
\sum_{n=-\infty}^{\infty} |c_n|^2 &=& \frac{1}{\pi} \int^{\pi}_{-\pi} |f(x)|^2  dx \nonumber \\
(\frac{a_0}{2})^2 + \sum_{n=-\infty}^{\infty} |a_n|^2 &=& \frac{1}{\pi} \int^{\pi}_{-\pi} x^4  dx \nonumber
\end{eqnarray}

Firstly, consider the left hand side of Parseval's Identity:
\begin{eqnarray}
(\frac{a_0}{2})^2 + \sum_{n=1}^{\infty} |a_n|^2 &=& \frac{2\pi^4}{9} + \sum_{n=1}^{\infty} \left[\frac{4(-1)^2n}{n^2}\right]^2 \nonumber \\
&=& \frac{2\pi^4}{9} + \sum_{n=1}^{\infty} \frac{16(-1)^2n}{n^4} \nonumber \\
&=& \frac{2\pi^4}{9} + 16 \sum_{n=1}^{\infty} \frac{1}{n^4} \nonumber \\
&=& \frac{2\pi^4}{9} + 16 \zeta (4) \nonumber
\end{eqnarray}

Secondly, consider the right hand side of Parseval's Identity:
\begin{eqnarray}
\frac{2}{\pi} \int^{\pi}_{0} x^4 dx &=& \frac{2}{\pi} \left[\frac{x^5}{5}\right]^{\pi}_{0} \nonumber \\
&=& \frac{2\pi^4}{5} \nonumber
\end{eqnarray}

Lastly, equate the left hand side to the right hand side of Parseval's Identity
\begin{eqnarray}
16\zeta (4) + \frac{2\pi^4}{9} &=& \frac{2\pi^4}{5} \nonumber \\
16\zeta (4) &=& \frac{8\pi^4}{45} \nonumber \\
\zeta (4) &=& \frac{\pi^4}{90} \nonumber
\end{eqnarray}

Now, using the definition of the zeta function we then see that

\begin{eqnarray}
\zeta (4) = \sum^{\infty}_{k=1} \frac{1}{k^4} = \frac{1}{1^4} + \frac{1}{2^4} +\frac{1}{3^4} + ... \nonumber
\end{eqnarray}

Multiplying by 90 and taking the 4th root we attain:

\begin{eqnarray}
{\pi^4} &=& 90\left[\frac{1}{1^4} + \frac{1}{2^4} +\frac{1}{3^4} + ... \right] \nonumber \\ \nonumber \\
\pi &=& \sqrt[4]{90\left[\frac{1}{1^4} + \frac{1}{2^4} +\frac{1}{3^4} + ... \right] \nonumber} \nonumber
\end{eqnarray}

\newpage
\subsection{Generating a value for $\zeta$ (6)}

Consider the odd function where $f(x) = x^3$ \\
Therefore, since $f(x)$ is an odd function then $a_n = 0$ where $\forall$ = 1,2,3 ... \nonumber \\

Now consider $b_n$,

\begin{eqnarray}
b_n = \frac{1}{\pi} \int^{\pi}_{-\pi} x^3sin(nx) dx &=& \frac{1}{\pi} \left[\frac{-x^3cos(nx)}{n}\right]^{\pi}_{-\pi} -\int^{\pi}_{-\pi} \frac{-3x^2cos(nx)}{n} dx \nonumber \\
&=& \frac{1}{\pi} \left[\left[\frac{-x^3cos(nx)}{n}\right]^{\pi}_{-\pi} + \left[\frac{3x^2sin(nx)}{n^2}\right]^{\pi}_{-\pi} - \int^{\pi}_{-\pi} \frac{6xsin(nx)}{n^2} dx \right] \nonumber \\ 
&=& \frac{1}{\pi} \left[\left[\frac{-x^3cos(nx)}{n}\right]^{\pi}_{-\pi} + \left[\frac{6xcos(nx)}{n^3}\right]^{\pi}_{-\pi} - \int^{\pi}_{-\pi} \frac{6cos(nx)}{n^3} dx \right] \nonumber \\
&=& \frac{1}{\pi} \left[\frac{-n^2x^3cos(nx)}{n^3} + \frac{6xcos(nx)}{n^3}\right]^{\pi}_{-\pi} \nonumber \\
&=& \frac{1}{\pi} \left[\frac{6xcos(nx) - n^2x^3cos(nx)}{n^3}\right]^{\pi}_{-\pi} \nonumber \\
&=& \frac{1}{\pi} \left[\frac{6\pi cos(n\pi) - \pi^3 n^2 cos(n\pi)}{n^3} - \frac{n^2\pi^3cos(-n\pi) - 6\pi cos(n\pi)}{n^3}\right] \nonumber \\
&=& \frac{1}{\pi} \left[\frac{6\pi cos(n\pi) - \pi^3 n^2 cos(n\pi) - n^2\pi^3cos(-n\pi) + 6\pi cos(n\pi))}{n^3} \right] \nonumber \\
&=& \frac{1}{\pi} \left[\frac{-2n^2\pi^3cos(n\pi) + 12\pi cos(n\pi)}{n^3} \right] \nonumber \\
&=& \left[\frac{12(-1)^n - 2\pi^2n^2(-1)^n}{n^3} \right] \nonumber \\
&=& \frac{12(-1)^n}{n^3} - \frac{2\pi^2 (-1)^n}{n} \nonumber
\end{eqnarray}

Now, by substituting $b_n$ back into the Fourier Series we obtain:

\begin{eqnarray}
x^3 = \sum^{\infty}_{n=1} \left[\frac{12(-1)^n}{n^3} - \frac{2\pi^2(-1)^n}{n}\right] sin(nx)  \mbox{     for     }    0 \leq x\leq \pi \nonumber
\end{eqnarray}

Then, consider Parsevals Identity where: \\
\begin{eqnarray}
c_n = a_n + ib_n \rightarrow |c_n|^2 &=& a_n^2 + b_n^2  \mbox{    (Note: $a_n = 0$ where $f(x)$ is odd)    } \nonumber \\
\sum_{n=-\infty}^{\infty} |c_n|^2 &=& \frac{1}{\pi} \int^{\pi}_{-\pi} |f(x)|^2  dx \nonumber \\
\sum_{n=-\infty}^{\infty} |b_n|^2 &=& \frac{1}{\pi} \int^{\pi}_{-\pi} x^6  dx \nonumber
\end{eqnarray}

Firstly, consider the left hand side of Parseval's Identity:
\begin{eqnarray}
\sum_{n=1}^{\infty} |b_n|^2 &=& \sum^{\infty}_{n=1} \left[\frac{12(-1)^n}{n^3} - \frac{2\pi^2(-1)^n}{n}\right]^2 \nonumber \\
&=& \sum_{n=1}^{\infty} \left[\frac{144(-1)^{2n}}{n^6} - \frac{48(-1)^{2n}}{n^4} + \frac{4\pi^4(-1)^{2n}}{n^2}\right] \nonumber \\
&=& 144 \sum_{n=1}^{\infty} \frac{1}{n^6} - 48\pi^2 \sum_{n=1}^{\infty} \frac{1}{n^4} + 4\pi^4 \sum_{n=1}^{\infty} \frac{1}{n^2} \nonumber \\
&=& 144 \zeta (6) - 48\pi^2 \zeta (4) + 4\pi^4 \zeta (2) \nonumber \\
&=& 144 \zeta(6) - \frac{8\pi^6}{15} + \frac{4\pi^6}{6} \nonumber
\end{eqnarray}

Secondly, consider the right hand side of Parseval's Identity:
\begin{eqnarray}
\frac{1}{\pi} \int^{\pi}_{-\pi} x^6 dx &=& \frac{1}{\pi} \left[\frac{x^7}{7}\right]^{\pi}_{-\pi} \nonumber \\
&=& \frac{2\pi^6}{7} \nonumber
\end{eqnarray}
  
Lastly, equate the left hand side to the right hand side of Parseval's Identity
\begin{eqnarray}
144\zeta (6) - \frac{8\pi^6}{15} + \frac{4\pi^6}{6}  &=& \frac{2\pi^6}{7} \nonumber \\
144\zeta (4) &=& \frac{16\pi^6}{105} \nonumber \\
\zeta (6) &=& \frac{\pi^6}{945} \nonumber
\end{eqnarray}

Now, using the definition of the zeta function we then see that

\begin{eqnarray}
\zeta (6) = \sum^{\infty}_{k=1} \frac{1}{k^6} = \frac{1}{1^6} + \frac{1}{2^6} +\frac{1}{3^6} + ... \nonumber
\end{eqnarray}

Multiplying by 945 and taking the 6th root we attain:

\begin{eqnarray}
{\pi^6} &=& 945\left[\frac{1}{1^6} + \frac{1}{2^6} +\frac{1}{3^6} + ... \right] \nonumber \\ \nonumber \\
\pi &=& \sqrt[6]{945\left[\frac{1}{1^6} + \frac{1}{2^6} +\frac{1}{3^6} + ... \right] \nonumber} \nonumber
\end{eqnarray}

\newpage
\subsection{Generating a value for $\zeta$ (8)}

Consider the even function where $f(x) = x^4$ \\
Therefore, since f(x) is an even function then $b_n = 0$ where $\forall$ = 1,2,3 ... \nonumber \\ 

First consider $a_0$,
\begin{eqnarray}
a_0 = \frac{1}{\pi} \int^{\pi}_{-\pi} x^4 dx = \frac{2}{\pi} \int^{\pi}_{0} x^4 dx = \frac{2\pi^4}{5} \nonumber
\end{eqnarray}

Then consider $a_n$,
\begin{eqnarray}
a_n = \frac{1}{\pi} \int^{\pi}_{-\pi} x^4 cos(nx) dx &=& \frac{2}{\pi} \int^{\pi}_{0} x^4 cos(nx) dx\nonumber \\
&=& \frac{2}{\pi} \left[\left[\frac{x^4 sin(nx)}{n}\right]^{\pi}_{0} - \int^{\pi}_{0} \frac{4x^3sin(nx)}{n} dx\right] \nonumber \\
&=& \frac{-2}{\pi} \int^{\pi}_{0} \frac{4x^3sin(nx)}{n}dx \nonumber \\
&=& \frac{-8}{n\pi} \int^{\pi}_{0} x^3 sin(nx) dx \nonumber \\
&=& \frac{-8}{n\pi} \left[\left[\frac{-x^3cos(nx)}{n}\right]^{\pi}_{0} + \int^{\pi}_{0} \frac{3x^2cos(nx)}{n} dx\right] \nonumber \\
&=& \frac{-8}{n\pi} \left[\left[\frac{-x^3cos(nx)}{n}\right]^{\pi}_{0} + \left[\frac{3x^2 sin(nx)}{n^2}\right]^{\pi}_{0} - \int^{\pi}_{0} \frac{6xsin(nx)}{n^2} dx \right] \nonumber \\
&=& \frac{-8}{n\pi} \left[\left[\frac{-x^3cos(nx)}{n}\right]^{\pi}_{0} + \left[\frac{6x cos(nx)}{n^3}\right]^{\pi}_{0} - \int^{\pi}_{0} \frac{6cos(nx)}{n^3} dx \right] \nonumber \\
&=& \frac{-8}{n\pi} \left[\frac{-n^2x^3cos(nx) + 6xcos(nx)}{n^3}\right]^{\pi}_{0} \nonumber \\
&=& \frac{-8}{n\pi} \left[\frac{-n^2\pi^3cos(n\pi) + 6\pi cos(n\pi)}{n^3}\right] \nonumber \\
&=& \left[\frac{8n^2\pi^2 cos(n\pi) - 48cos(n\pi)}{n^4}\right] \nonumber \\
&=& \frac{8\pi^2(-1)^n}{n^2} - \frac{48(-1)^n}{n^4} \nonumber
\end{eqnarray}

Now, by substituting back $a_n$ and $a_0$ back into the Fourier Series we obain:

\begin{eqnarray}
x^4 = \frac{\pi^4}{5} + \sum_{n=1}^{\infty} \left[\frac{8n^2\pi^2(-1)^n}{n^4} -\frac{48(-1)^n}{n^4}\right] cos(nx)     \mbox{     for     }    -\pi\leq x\leq \pi \nonumber
\end{eqnarray}

Then, consider Parsevals Identity where: \\
\begin{eqnarray}
c_n = a_n + ib_n \rightarrow |c_n|^2 &=& a_n^2 + b_n^2  \mbox{    (Note: $b_n = 0$ where $f(x)$ is even)    } \nonumber \\
\sum_{n=-\infty}^{\infty} |c_n|^2 &=& \frac{1}{\pi} \int^{\pi}_{-\pi} |f(x)|^2  dx \nonumber \\
(\frac{a_0}{2})^2 +\sum_{n=-\infty}^{\infty} |a_n|^2 &=& \frac{1}{\pi} \int^{\pi}_{-\pi} x^8  dx \nonumber
\end{eqnarray}

Firstly, consider the left hand side of Parseval's Identity:
\begin{eqnarray}
(\frac{a_0}{2})^2 + \sum_{n=1}^{\infty} |a_n|^2 &=& \frac{2\pi^8}{25} + \sum_{n=1}^{\infty} \left[\frac{8\pi^2 (-1)^ n}{n^2} - \frac{48(-1)^n}{n^4}\right]^2 \nonumber \\
&=& \frac{2\pi^8}{25} + \sum_{n=1}^{\infty} \frac{64\pi^4(-1)^{2n}}{n^4} - \frac{784\pi^2(-1)^{2n}}{n^6} + \frac{2304(-1)^{2n}}{n^8} \nonumber \\
&=& \frac{2\pi^8}{25} + 64\pi^4 \sum_{n=1}^{\infty} \frac{1}{n^4} - 784\pi^2 \sum_{n=1}^{\infty} \frac{1}{n^6} + 2304\sum_{n=1}^{\infty} \frac{1}{n^8} \nonumber \\
&=& \frac{2\pi^8}{25} + 64\pi^4 \zeta (4) - 784\pi^2 \zeta (6) + 2304 \zeta (8) \nonumber \\
&=& \frac{2\pi^8}{25} + \frac{64\pi^8}{90} - \frac{768\pi^8}{945} + 2304 \zeta (8) \nonumber
\end{eqnarray}

Secondly, consider the right hand side of Parsevals Identity:
\begin{eqnarray}
\frac{2}{\pi} \int^{\pi}_{0} x^8 dx &=& \frac{2}{\pi} \left[\frac{x^9}{9}\right] \nonumber \\
&=& \frac{2\pi^8}{9} \nonumber
\end{eqnarray}

Lastly, equate the left hand side to the right hand side of Parseval's Identity
\begin{eqnarray}
\frac{2\pi^8}{25} + \frac{64\pi^8}{90} - \frac{768\pi^8}{945} + 2304 \zeta (8) &=& \frac{2\pi^8}{9} \nonumber \\
2304 \zeta (8) &=& \frac{128\pi^8}{525}  \nonumber \\
\zeta (8) &=& \frac{\pi^8}{9450} \nonumber
\end{eqnarray}

Now, using the definition of the zeta function we then see that

\begin{eqnarray}
\zeta (8) = \sum^{\infty}_{k=1} \frac{1}{k^8} = \frac{1}{1^8} + \frac{1}{2^8} +\frac{1}{3^8} + ... \nonumber
\end{eqnarray}

Multiplying by 9450 and taking the 8th root we attain:

\begin{eqnarray}
{\pi^8} &=& 9450\left[\frac{1}{1^8} + \frac{1}{2^8} +\frac{1}{3^8} + ... \right] \nonumber \\ \nonumber \\
\pi &=& \sqrt[8]{9450\left[\frac{1}{1^8} + \frac{1}{2^8} +\frac{1}{3^8} + ... \right] \nonumber} \nonumber
\end{eqnarray}

\newpage
\section{The Wallis Product}

Consider the function in which $f(x) = sin^n (x)$

\begin{eqnarray}
I(n) = \int^{\pi}_{0} sin^n (x) dx \nonumber
\end{eqnarray}

By evaluating this integral we obtain the following:

\begin{eqnarray}
I(n) = \int^{\pi}_{0} sin^n (x) dx &=& \int^{\pi}_{0} sin^{n-1}(x)sin(x) dx \nonumber \\
&=& \left[sin^{n-1}(x)(-cos(x))\right]^{\pi}_{0} - \int^{\pi}_{0} (n-1)sin^{n-2}(x)(-cos(x))(cos(x))dx \nonumber \\
&=& 0 + (n-1)\int^{\pi}_{0} sin^{n-2}(x)cos^(x) dx \nonumber \\
&=& (n-1)\int^{\pi}_{0}sin^{n-2}x (1-sin^{2}x) dx\nonumber \\
&=& (n-1)[I(n-2)-I(n)] \nonumber
\end{eqnarray}

Now it is possible to generate values that are either even or odd.
To generate even values, consider the following:

\begin{eqnarray}
I(n) &=& (n-1)[I(n-2)-I(n)] \nonumber \\
I(n) &=& (n-1)I(n-2)-(n-1)I(n) \nonumber \\
I(n)+(n-1)I(n) &=& (n-1)I(n-2) \nonumber \\
I(n)[1 + (n-1)] &=& (n-1)I(n-2) \nonumber \\
I(n) &=& \frac{n-1}{n} I(n-2) \mbox{                 [1]} \nonumber
\end{eqnarray}

 To generate odd values, consider the following: By letting $n=2m+1$ in $[1]$, we attain the following:

\begin{eqnarray}
I(2m+1) &=& \frac{2m}{2m+1} I(2m-1) \mbox{                [2]} \nonumber \\
\frac{I(2m+1)}{I(2m-1)} &=& \frac{2m}{2m+1} \nonumber \\
\frac{I(2m-1)}{I(2m+1)} &=& \frac{2m+1}{2m} \mbox{                [3]} \nonumber
\end{eqnarray}

\newpage
Now we can consider the product of values for both even and odd functions

\begin{eqnarray}
I(0) = \int^{\pi}_{0} sin^0 (x) dx = \int^{\pi}_{0} dx = \pi \nonumber
\end{eqnarray}

\begin{eqnarray}
I(1) = \int^{\pi}_{0} sin(x) dx = [-cos(x)]^{\pi}_{0} \int^{\pi}_{0} dx = -(-1)-(-1) = 2 \nonumber
\end{eqnarray}

Note: for I(2n), $n=2m$ is substitued into [1]

\begin{eqnarray}
I(2m) = \int^{\pi}_{0} sin^{2m}(x) dx = \frac{2m-1}{2m} . \frac{2m-3}{2m-2} . I(2m-4)\nonumber
\end{eqnarray}

Note: for I(2n+1), $n=2m+1$ is substituted into [2]

\begin{eqnarray}
I(2m+1) = \int^{\pi}_{0} sin^{2m+1}(x) dx = \frac{2m}{2m+1} . I(2m-1)\nonumber
\end{eqnarray}

For even functions, consider I(2n) and expand further by repeating the iterative process.

\begin{eqnarray}
I(2n) = \frac{2m-1}{2m}.\frac{2m-3}{2m-2}.\frac{2m-5}{2m-4} ... \frac{5}{6}.\frac{3}{4}.\frac{1}{2}.I(0) = \pi \prod^{m}_{k=1} \frac{2k-1}{2k} [4]\nonumber
\end{eqnarray}

For odd functions, consider I(2n+1) and expand further by repeating the iterative process.

\begin{eqnarray}
I(2n+1) = \frac{2m}{2m+1}.\frac{2m-2}{2m-1}.\frac{2m-1}{2m-3}.\frac{6}{7}.\frac{4}{5}.\frac{2}{3}.I(1) = 2 \prod^{m}_{k=1} \frac{2k}{2k+1}  [5] \nonumber
\end{eqnarray}

Therefore,

\begin{eqnarray}
\frac{I(2n)}{I(2n+1)} = \frac{\pi}{2} \prod^{\infty}_{k=1} \frac{2k-1}{2k}.\frac{2k+1}{2k} \nonumber
\end{eqnarray}

From this, we can deduce the following as $0\leq sinx \leq 1$:

\begin{eqnarray}
sin^{2n+1}(x)]\leq sin^{2n}(x)\leq sin^{2n-1}(x) \mbox{     for     } 0 \leq x \leq \pi \nonumber \\ \int^{\pi}_{0}sin^{2n+1}(x) dx \leq \int^{\pi}_{0}sin^{2n}(x) dx \leq \int^{\pi}_{0}sin^{2n-1}(x) dx \nonumber \\
I(2n+1)\leq I(2n) \leq I(2n-1) \nonumber \\
I \leq \frac{I(2n)}{I(2n+1)} \leq \frac{I(2n-1)}{I(2n+1)} \nonumber
\end{eqnarray}

Then, by substituting $[3])$ for $$\frac{I(2n-1)}{I(2n+1)}$$ 
we obtain:

\begin{eqnarray}
1\leq \frac{I(2n)}{I(2n+1)} \leq \frac{2n+1}{2n} \nonumber
\end{eqnarray} \\

Now consider the Squeeze Theorem which states the following:

\begin{eqnarray}
\mbox{If $a_n$, $b_n$ and $c_n$ exist for all $a_n \leq b_n \leq c_n$ for all n} \nonumber \\
\lim_{n\to\infty} a_n = \lim_{n\to\infty} = \lim_{n\to\infty} c_n &=& L < \infty \nonumber \\
\mbox{   Then   } \lim_{n\to\infty} b_n &=& L \nonumber
\end{eqnarray}

In this case, the lower limit would be the following:

\begin{eqnarray} 
\lim_{n\to\infty} a_n = 1 \nonumber  
\end{eqnarray}

Now, consider the limit of $c_n$,

\begin{eqnarray}
\lim_{n\to\infty} \frac{I(2n-1)}{I(2n+1)} = \lim_{n\to\infty} \frac{2n+1}{2n} = \lim_{n\to\infty} \left(1 + \frac{1}{2n}\right)  \rightarrow 1 \nonumber
\end{eqnarray}

Henceforth, since both $a_n$ and $c_n$ tend towards the same limit which in this case is 1 then $b_n$ must also tend towards 1 by the squeeze theorem which is shown below.

\begin{eqnarray}
\lim_{n\to\infty} \left[\frac{I(2n)}{I(2n+1)}\right] &=& 1 \nonumber \\
\frac{\pi}{2} \lim_{n\to\infty} \prod^{n}_{k=1} \frac{2k-1}{2k} . \frac{2k+1}{2k} &=& 1 \nonumber
\end{eqnarray}

Therefore,
\begin{eqnarray}
\frac{\pi}{2} &=& \prod^{\infty}_{k=1} \left[\frac{2k}{2k-1} . \frac{2k}{2k+1}\right] = \frac {2}{1}.\frac{2}{3}.\frac{4}{3}.\frac{4}{5}.\frac{6}{5}.\frac{6}{7} ... \nonumber \\
\pi &=& \prod^{\infty}_{k=1} 2\left[\frac{2k}{2k-1} . \frac{2k}{2k+1}\right] = 2[\frac {2}{1}.\frac{2}{3}.\frac{4}{3}.\frac{4}{5}.\frac{6}{5}.\frac{6}{7} ...] \nonumber
\end{eqnarray}

And lastly,

\begin{eqnarray}
\pi = 2[\frac{2}{1}.\frac{2}{3}.\frac{4}{3}.\frac{4}{5}.\frac{6}{5}.\frac{6}{7}.\frac{8}{7}.\frac{8}{9} ...] \cite{cite_key12} \nonumber
\end{eqnarray}

It is also important to acknowledge that this method to calculate $\pi$ was different to the other methods used to calculate $\pi$ because it is a different type of infinite decent. Unlike Newtons Arcsine function with the usage of an infinite sum, The Wallis product is a unique infinite product that convereges towards $\pi$ which was made without Infinitesimal calculus as it did not exist at the time.

\newpage
\section{Continued Fractions}
\subsection{Preface list of data and proof's used in the calculations}
\subsubsection{The Taylor Series}

Provided that $f(x)$ is infinitely differentiable at x=a,
\begin{eqnarray}
f(x)f(a) + \frac{f'(a)}{1!}(x-a) + \frac{f''(a)}{2!}(x-a)^2 + \frac{f'''(a)}{3!}(x-a)^3 + ... \nonumber
\end{eqnarray}

\subsubsection{Euler's Continued Fraction Formula}

\begin{eqnarray}
a_0 + a_0a_1 + a_0a_1a_2 + ... + a_0a_1a_2...a_n = 
\cfrac{a_0}{\displaystyle 1-
\cfrac{a_1}{\displaystyle {1+a_1} -
\cfrac{a_2}{\displaystyle {1+a_2} - 
\cfrac{a_3\ddots}{\displaystyle {\ddots1+a_{n-1}} -
\cfrac{a_n}{1+a_n}}}}} \nonumber
\end{eqnarray}

\subsubsection{The Natural Logarithm in Continued Fractions}

Let $f(z)=ln\frac{1+z}{1-z}$ and then consider the Taylor Expansion of the function to develop an infinite series

\begin{eqnarray}
f'(z) &=& \frac{-2}{z^2 - 1} \nonumber \\
f''(z) &=& \frac{4z}{(1-z^2)^2} \nonumber \\
f'''(z) &=& \frac{16z^2}{(1-z^2)^3} + \frac{4}{(1-z^2)^2} \nonumber \\
f^{(4)}(z) &=& \frac{96z^3}{(1-z^2)^4} +\frac{48z}{(1-z^2)^3} \nonumber \\
f^{(5)}(z) &=& \frac{768z^4}{(1-z^2)^5} + \frac{576z^2}{(1-z^2)^4} + \frac{48}{(1-z^2)^3} \nonumber
\end{eqnarray}

Then, by letting $a = 0$ in and expanding through the Taylor Series the following is obtained:

\begin{eqnarray}
ln\frac{1+z}{1-z} &=& \frac{2z}{1!} + \frac{4z^3}{3!} + \frac{48z^5}{5!} ... \nonumber \\
&=& 2[z + \frac{z^3}{3} + \frac{z^5}{5} + \frac{z^7}{7} ... ] \nonumber \\
&=& 2 \sum_{n=0}^{\infty} \frac{z^{2n+1}}{2n+1} \nonumber
\end{eqnarray}

By extracting the series formed in the previous calculation:

\begin{eqnarray}
ln\frac{1+z}{1-z} = 2z[1 + \frac{z^2}{3} + \frac{z^2}{3}\frac{z^2}{(5/3)} + \frac{z^2}{3}\frac{z^2}{(5/3)}\frac{z^2}{(7/5)} + \frac{z^2}{3}\frac{z^2}{(5/3)}\frac{z^2}{(7/5)}\frac{z^2}{(9/7)} + ... \nonumber
\end{eqnarray}

The following Continued Fraction can be produced by using Eulers Formula

\begin{eqnarray}
\cfrac{2z}{1-\cfrac{\frac{z^2}{3}}{1+\frac{z^3}{3}-\cfrac{\frac{3z^2}{5}}{1+\frac{3z^2}{5}-\cfrac{\frac{5z^2}{7}}{1+\frac{5z^2}{7}-\cfrac{\frac{7z^2}{9}}{1+\frac{7z^2}{9}-\cfrac{9z^2}{11}}}}}} \nonumber
\end{eqnarray}

\subsection{A Continued Fraction for $\pi$}

Consider the complex number where $z = \frac{1+i}{1-i}$ = i (By muliplying by the conjugate)

Now, consider the complex number in Euler Form,

\begin{eqnarray}
\frac{1+i}{1-i} &=& e^{\frac{i\pi}{2}} \nonumber \\
ln\frac{1+i}{1-i} &=& \frac{i\pi}{2} \nonumber \\
\pi &=& \frac{2}{i} ln\frac{1+i}{1-i} \nonumber
\end{eqnarray}

Letting z = i in Eulers Continued Fraction formula for the Natural Logarithm we can attain the following continued fraction

\begin{eqnarray}
ln\frac{1+i}{1-i} = \cfrac{2i}{1+\cfrac{\frac{1}{3}}{1-\frac{1}{3}+\cfrac{\frac{3}{5}}{1-\frac{3}{5}+\cfrac{\frac{5}{7}}{1-\frac{5}{7}+\cfrac{\frac{7}{9}}{1-\frac{7}{9}+\cfrac{\frac{9}{11}}{\ddots}}}}}} \nonumber
\end{eqnarray}

\begin{eqnarray}
ln\frac{1+i}{1-i} = \cfrac{2i}{1+\cfrac{1}{3-1+\cfrac{\frac{3^2}{5}}{1-\frac{3}{5}+\cfrac{\frac{5}{7}}{1-\frac{5}{7}+\cfrac{\frac{7}{9}}{1-\frac{7}{9}+\cfrac{\frac{9}{11}}{\ddots}}}}}} \nonumber
\end{eqnarray}

\begin{eqnarray}
ln\frac{1+i}{1-i} = \cfrac{2i}{1+\cfrac{1}{2+\cfrac{3^2}{5-3\cfrac{\frac{5^2}{7}}{1-\frac{5}{7}+\cfrac{\frac{7}{9}}{1-\frac{7}{9}+\cfrac{\frac{9}{11}}{\ddots}}}}}} \nonumber
\end{eqnarray}

And by repeating this process the following is eventually attained:

\begin{eqnarray}
ln\frac{1+i}{1-i} = \cfrac{2i}{1+\cfrac{1}{2+\cfrac{3^2}{2+\cfrac{5^2}{2+\cfrac{7^2}{2+\cfrac{9^2}{\ddots}}}}}} \nonumber
\end{eqnarray}

Now, by multiplying by $(\frac{i}{2})$, the following is attained:

\begin{eqnarray}
\pi = (\frac{i}{2})ln\frac{1+i}{1-i} = \cfrac{4}{1^2+\cfrac{1^2}{2+\cfrac{3^2}{2+\cfrac{5^2}{2+\cfrac{7^2}{2+\cfrac{9^2}{\ddots}}}}}} \cite{cite_key8} \nonumber
\end{eqnarray}

\newpage
\section{The Gregory-Leibniz Series}
\subsection{Preface list of data and proof's used in the calculations}
\subsubsection{The Taylor Series}

Provided that $f(x)$ is infinitely differentiable at x=a,
\begin{eqnarray}
f(x)f(a) + \frac{f'(a)}{1!}(x-a) + \frac{f''(a)}{2!}(x-a)^2 + \frac{f'''(a)}{3!}(x-a)^3 + ... \nonumber
\end{eqnarray}

\subsection{$Arctan(x)$}

Consider the function $f(x) = arctan(x)$,
By using the Taylor expansion, we can develop this function into an infinite series.

\begin{eqnarray}
f(x) &=& Arctan(x) \nonumber \\
f'(x) &=& \frac{1}{(1+x^2)} \nonumber \\
f''(x) &=& \frac{-2x}{(1+x^2)^2} \nonumber \\
f'''(x) &=& \frac{8x^2}{(1+x^2)^3} - \frac{2}{(1+x^2)^2} \nonumber \\
f^{4}(x) &=& \frac{-48x^3}{(1+x^2)^4} + \frac{24x}{(1+x^2)^3} \nonumber \\
f^{5}(x) &=& \frac{384x^4}{(1+x^2)^5} - \frac{288x^2}{(1+x^2)^4} + \frac{24x}{(1+x^2)^3} \nonumber \\
f^{6}(x) &=& \frac{-3840^6}{(1+x^2)^6} + \frac{3840x^3}{(1+x^2)^5} - \frac{720x}{(1+x^2)^4} \nonumber \\
f^{7}(x) &=& \frac{46080^6}{(1+x^2)^7} - \frac{57600x^4}{(1+x^2)^6} + \frac{17280x^2}{(1+x^2)^5} - \frac{720}{(1+x^2)^4} \nonumber
\end{eqnarray}

With further expansion through the Taylor Series with $a=0$, the following can be obtained:

\begin{eqnarray}
Arctan(x) &=& x - \frac{2!x^3}{3!} + \frac{4!x^5}{5!} - \frac{6!x^7}{7!} ... \nonumber \\
Arctan(x) &=& x - \frac{x^3}{3} + \frac{x^5}{5} - \frac{x^7}{7} ... \nonumber \\
Arctan(x) &=& \sum_{n=0}^{\infty} \frac{(-1)^n}{2n+1}x^{2n+1} \cite{cite_key2} \nonumber
\end{eqnarray}

Now, consider when $x = 1$.

\begin{eqnarray}
\frac{\pi}{4} &=& 1 - \frac{1}{3} +\frac{1}{5} - \frac{1}{7} + \frac{1}{9} - \frac{1}{11} + \frac{1}{13} - \frac{1}{15} + \frac{1}{17} ... \nonumber \\
\pi &=& 4(1 - \frac{1}{3} +\frac{1}{5} - \frac{1}{7} + \frac{1}{9} - \frac{1}{11} + \frac{1}{13} - \frac{1}{15} + \frac{1}{17} ...) \cite{cite_key2} \nonumber
\end{eqnarray}

\newpage
\section{The Newtons Series expansion of the Arcsine Function}
\subsection{Preface list of data and proof's used in the calculations}
\subsubsection{The Taylor Series}

Provided that $f(x)$ is infinitely differentiable at x=a,
\begin{eqnarray}
f(x)f(a) + \frac{f'(a)}{1!}(x-a) + \frac{f''(a)}{2!}(x-a)^2 + \frac{f'''(a)}{3!}(x-a)^3 + ... \nonumber
\end{eqnarray}

\subsection{$Arcsine(x)$}

Consider the function $f(x) = arctan(x)$,
By using the Taylor expansion, we can develop this function into an infinite series

\begin{eqnarray}
f(x) &=& Arcsine(x) \nonumber \\
f'(x) &=& \frac{1}{\sqrt{1-x^2}} \nonumber \\
f''(x) &=& \frac{x}{(1-x^2)^\frac{3}{2}} \nonumber \\
f'''(x) &=& \frac{3x^2}{(1-x^2)^\frac{5}{2}} + \frac{1}{(1-x^2)^\frac{3}{2}} \nonumber \\
f^{4}(x) &=& \frac{15x^3}{(1-x^2)^\frac{7}{2}} + \frac{9x}{(1-x^2)^\frac{5}{2}} \nonumber \\
f^{5}(x) &=& \frac{105x^4}{(1-x^2)^\frac{9}{2}} + \frac{90x^2}{(1-x^2)^\frac{7}{2}} + \frac{9}{(1-x^2)^\frac{5}{2}} \nonumber \\
f^{6}(x) &=& \frac{945^5}{(1-x^2)^\frac{11}{2}} + \frac{1050x^3}{(1-x^2)^\frac{9}{2}} - \frac{720x}{(1+x^2)^\frac{5}{2}} \nonumber \\
f^{7}(x) &=& \frac{10395^6}{(1-x^2)^\frac{13}{2}} + \frac{14175x^4}{(1-x^2)^\frac{11}{2}} + \frac{4725x^2}{(1-x^2)^\frac{9}{2}} + \frac{225}{(1-x^2)^\frac{7}{2}} \nonumber
\end{eqnarray}

With further expansion through the Taylor Series with $a=0$, the following can be obtained:

\begin{eqnarray}
Arcsine(x) &=& \frac{x}{1!} + \frac{x^3}{3!} + \frac{9x^5}{5!} + \frac{225x^7}{7!} + ...\nonumber \\
Arcsine(x) &=& x + \frac{x^3}{6} + \frac{3x^5}{40} + \frac{5x^7}{112}\nonumber + ... \\
Arcsine(x) &=& \sum_{n=0}^{\infty} \frac{(2n)!}{2^{2n}(n!)^2(2n+1)}x^{2n+1} \cite{cite_key3} \nonumber
\end{eqnarray}

Now, consider when $x= \frac{1}{2}$

\begin{eqnarray}
\frac{\pi}{6} &=& \frac{1}{2} + \frac{1}{48} + \frac{3}{40} + \frac{15}{43008} + ... \nonumber \\ 
\pi &=& 6(\frac{1}{2} + \frac{1}{48} + \frac{3}{40} + \frac{15}{43008} + ...) \cite{cite_key3} \nonumber
\end{eqnarray}

\newpage
\section{Viete's Formula with Nested Radicals}

First, Start by recalling the simple trigonometric identity in which $sin(2x)=2sin(x)cos(x)$ \\

Now consider $sin(2^n x)$,

\begin{eqnarray}
sin(2^n x) &=& 2cos(2^{n-1}x)sin(2^{n-1}x) \nonumber \\
&=& 2cos(2^{n-1}x)[2cos(2^{n-2}x)sin(2^{n-2}x) ...] \nonumber \\ \nonumber \\ 
sin(2^n x) &=& 2^n [cos(2^{n-1}x)cos(2^{n-2}x) ... cos(2x)cos(x)]sin(x) \nonumber \\
\frac{sin(2^n x)}{2^n sin(x)} &=& cos(2^{n-1}x)cos(2^{n-2}x) ... cos(2x)cos(x) \nonumber \\
&=& \prod^{n-1}_{i=0} cos(2^i x) \nonumber
\end{eqnarray}

Now, let $x = \frac{y}{2^n}$

\begin{eqnarray}
\frac{sin(\frac{2^n y}{2^n})}{2^n sin(\frac{y}{2^n})} &=& \prod^{n-1}_{i=0} cos(\frac{2^i y}{2^n}) \nonumber \\
\frac{sin(y)}{2^n sin(\frac{y}{2^n})} &=& \prod^{n-1}_{i=0} cos(\frac{2^i y}{2^n j}) \nonumber \\
&=& \prod^{n}_{j=1} cos(\frac{y}{2^j }) \nonumber \\
&=& cos(\frac{y}{2})cos(\frac{y}{2^2})...cos(\frac{y}{2^n}) \nonumber
\end{eqnarray}

Therefore,

\begin{eqnarray}
\frac{2sin(\frac{y}{2})cos(\frac{y}{2})}{2^nsin(\frac{y}{2^n})cos(\frac{y}{2})} &=& cos(\frac{y}{2^2})...cos(\frac{y}{2^n}) \nonumber \\
\frac{2sin(\frac{y}{2})}{2^nsin(\frac{y}{2^n})} &=& \prod^{n}_{i=2} cos(\frac{\pi}{2^i}) \nonumber \\
\frac{2}{2^nsin(\frac{\pi}{2^n})} &=& \prod^{n}_{i=2} cos(\frac{\pi}{2^i}) \nonumber \\
\frac{\frac{2}{\pi}}{\frac{2^n}{\pi}sin(\frac{\pi}{2^n})} &=& \prod^{n}_{i=2} cos(\frac{\pi}{2^i}) \nonumber
\end{eqnarray}

\newpage

Now,  let $\lim_{n \rightarrow \infty} $for the above and obtain the following:

\begin{eqnarray}
\frac{2}{\pi} &=& \prod^{\infty}_{i=2} cos(\frac{\pi}{2^i}) \nonumber \\
\pi &=& 2\prod^{\infty}_{i=2} \frac{1}{cos(\frac{\pi}{2^i})} \nonumber
\end{eqnarray}

Now consider the following, by re-arranging the half cosine formula, we obtain an identity for $cos(\frac{\theta}{2})$

\begin{eqnarray}
cos(2\theta) &=& 2cos^{2} - 1 \nonumber \\
cos(\theta) &=& 2cos^{2}(\frac{\theta}{2}) - 1 \nonumber \\
2cos(\frac{\theta}{2}) &=& \sqrt{\frac{1+cos{\theta}}{2}} \nonumber \\
2cos(\frac{\theta}{2}) &=& \sqrt{2 + 2cos\theta} \nonumber \\
cos(\frac{\theta}{2}) &=& \frac{\sqrt{2 + 2cos\theta}}{2} \nonumber 
\end{eqnarray}

Therefore,

\begin{eqnarray}
cos(\frac{\pi}{2^2}) &=& \frac{\sqrt{2+2.0}}{2} = \frac{\sqrt2}{2} \nonumber \\
cos(\frac{\pi}{2^3}) &=& \frac{\sqrt{2+\sqrt2}}{2} \nonumber \\
cos(\frac{\pi}{2^4}) &=& \frac{\sqrt{2+{\sqrt{2+\sqrt2}}}}{2} \nonumber
\end{eqnarray}

Now by using this in the originial product after applying $\lim_{n \rightarrow \infty}$, the following is obtained:

\begin{eqnarray}
\prod^{n}_{i=2} cos(\frac{\pi}{2^i}) &=& \overbrace{\frac{\sqrt2}{2} \times \frac{\sqrt{2+\sqrt2}}{2} \times \frac{\sqrt{2+{\sqrt{2+\sqrt2}}}}{2}}^{\mbox{n-1 terms}} \times ... \nonumber \\
&=& \frac{\sqrt2 \times \sqrt{2+\sqrt2} \times \sqrt{2+{\sqrt{2+\sqrt2}}} \times ...}{2^{n-1}} \nonumber
\end{eqnarray}

\newpage

Therefore,

\begin{eqnarray}
\prod^{\infty}_{i=2} \frac{1}{cos(\frac{\pi}{2^i})} &=& \frac{2^{n-1}}{\sqrt2 \times \sqrt{2+\sqrt2} \times \sqrt{2+{\sqrt{2+\sqrt2}}}} \times ... \nonumber \\
2\prod^{\infty}_{i=2} \frac{1}{cos(\frac{\pi}{2^i})} &=& \underbrace{\frac{2^n}{\sqrt2 \times \sqrt{2+\sqrt2} \times \sqrt{2+{\sqrt{2+\sqrt2}}}}}_{\mbox{n-1 terms}} \times ... \nonumber
\end{eqnarray}

Using n=3 and n=4 respectivly as examples,

\begin{eqnarray}
\frac{2^3}{\sqrt2 \sqrt{2+\sqrt2}} = \frac{2^3 \sqrt{2-\sqrt2}}{\sqrt2 \sqrt{4-2}} = \frac{2^3 \sqrt{2-\sqrt2}}{2} = 2^2 \sqrt{2-\sqrt2} \nonumber
\end{eqnarray}

\begin{eqnarray}
\frac{2^4}{\sqrt2 \sqrt{2+{\sqrt2}} \sqrt{2+\sqrt{2+\sqrt2}}} &=& \frac{2^4 \sqrt{2-\sqrt{2+\sqrt2}}}{\sqrt2 \sqrt{2+{\sqrt2}} \sqrt{4- (2+\sqrt2)}}  \nonumber \\
&=& \frac{2^4 \sqrt{2-\sqrt{2+\sqrt2}}}{\sqrt2 \sqrt{2+{\sqrt2}} \sqrt{2-\sqrt2}}  \nonumber \\
&=& \frac{2^4 \sqrt{2-\sqrt{2+\sqrt2}}}{\sqrt2 \sqrt{4-2}}  \nonumber \\
&=& 2^{3} \sqrt{2-\sqrt{2+\sqrt2}} \nonumber
\end{eqnarray}

So,

\begin{eqnarray}
\underbrace{\frac{2^n}{\sqrt2 \times \sqrt{2+\sqrt2} \times \sqrt{2+{\sqrt{2+\sqrt2}}}\times ...}}_{\mbox{n-1 terms}} \nonumber &=& 2^{n-1} {\sqrt{2-\underbrace{\sqrt{2+\sqrt{2 + ...}}}_{\mbox{n-2 terms}}}} \\
\prod^{\infty}_{i=2} \frac{1}{cos(\frac{\pi}{2^i})} &=& 2^{n-1} {\sqrt{2-\underbrace{\sqrt{2+\sqrt{2 + ...}}}_{\mbox{n-2 terms}}}} \nonumber \\
\Rightarrow 2\prod^{n+2}_{i=2} \frac{1}{cos(\frac{\pi}{2^i})} &=& 2^{n+1} {\sqrt{2-\underbrace{\sqrt{2+\sqrt{2 + ...}}}_{\mbox{n terms}}}} \nonumber
\end{eqnarray}

Therefore,

\begin{eqnarray}
\pi =\lim_{n\to\infty} 2^{n+1}\sqrt{2-\underbrace{\sqrt{2+\sqrt{2 + ...}}}_{\mbox{n terms}}} \nonumber \cite{cite_key1}
\end{eqnarray}

\newpage
\section{Analysis with Data Tables and Graphs of the various methods}
\subsection{The Wallis Product}
\begin{table}[ht]
\caption{As $n$ increases how does the Wallis Product converege towards $\pi$}

\centering
\begin{tabular}{c c c}
\hline\hline 
n value & Method\#1 & Error (\%) \\ [.5ex]
\hline
5 & 3.002175954556907 & 4.43777\\
10 & 3.067703806643499 & 2.35196\\
15  & 3.091336888596228 & 1.59969\\
20  & 3.103516961539234 & 1.21199\\
25  & 3.110945166901554 & 0.97554\\
30  & 3.115948285887959 & 0.81629\\
35  & 3.119547206305518 & 0.70173\\ 
40  & 3.122260326421437 & 0.61537\\
45  & 3.124378835915516 & 0.54793\\
50  & 3.126078900215411 & 0.49382\\
55  & 3.127473350412857 & 0.44943\\
60  & 3.128637797891591 & 0.41237\\
65  & 3.129624812079802 & 0.38095\\
70  & 3.130472076319065 & 0.35398\\
75  & 3.131207308587379 & 0.33058\\
80  & 3.131851351372613 & 0.31008\\
85  & 3.132420179022906 & 0.29197\\
90  & 3.132926240627509 & 0.27586\\
95  & 3.133379381619937 & 0.26144\\
100  & 3.133787490628162 & 0.24845\\
1000  & 3.140807746030395 & 0.02498\\
10000 & 3.141514118681922 & 0.00250 \\
100000 & 3.141584799657247 & 0.00025\\
1000000  & 3.141591868192124 & 0.00008\\
10000000  & 3.141592575049982 & 0.00001\\ [1ex]
\hline
\end{tabular}
\end{table}

The following data tables have been generated using excel and the calculations have been done in code using wolfram mathematica. There is very high accuracy in these values because mathematica has the ability to calculate with to 31 decimal places however the accuracy of the attained values for the percentage errors decrease as the values of n increase because as each method converges towards $\pi$ the percentage error becomes more susceptible to rounding errors and at large n values the percentage error becomes so small that Mathematica has trouble calculating and a rounding error is produced as the percentage error is only taken to 5 decimal places. However, this error is insignificant and therefore will not effect the results drastically.

\newpage

\subsection{The Gregory-Leibniz Series}

\begin{table}[ht]
\caption{As $n$ increases how does the Gregory-Leibniz Series converge towards $\pi$}
\centering
\begin{tabular}{c c c}
\hline\hline 
n value & Method\#2 & Error(\%) \\ [0.5ex]
\hline
5 & 2.976046176046176 & 5.26951 \\
10 & 3.232315809405593 & 2.88781 \\
15 & 3.079153394197426 & 1.98750 \\
20 & 3.189184782277595 & 1.51490 \\
25 & 3.103145312886011 & 1.22382 \\
30 & 3.173842337190749 & 1.02654 \\
35 & 3.113820229023573 & 0.88402 \\
40 & 3.165979272843215 & 0.77625 \\
45 & 3.119856090062712 & 0.69110 \\
50 & 3.161198612987056 & 0.62408 \\
55 & 3.123736933726277 & 0.56837 \\
60 & 3.157984995168666 & 0.52178 \\
65 & 3.126442007766234 & 0.48226 \\
70 & 3.155676462307475 & 0.44830 \\
75 & 3.128435328236984 & 0.41881 \\
80 & 3.153937862272616 & 0.37012 \\
85 & 3.129965139593801 & 0.34978 \\
90 & 3.152581332875124 & 0.33156 \\
95 & 3.131176269454981 & 0.31515 \\
100 & 3.151493401070910 & 0.31118 \\
1000 & 3.142591654339543 & 0.00318\\
10000 & 3.141692643590543 & 0.00032 \\
100000 & 3.141602653489794 & 0.00003\\
1000000  & 3.141591868192127 & 0.00001\\
10000000  & 3.141592653518277 & 0.00000\\ [1ex]  
\hline
\end{tabular}
\end{table}

\newpage

\subsection{Newtons Arcsine Expansion Series for $\pi$}
\begin{table}[ht]
\caption{As $n$ increases how does Newtons Arcsine Expansion converge towards $\pi$}
\centering
\begin{tabular}{c c c}
\hline\hline 
n value & Method\#3 & Error(\%) \\ [0.5ex]
\hline
5 & 3.141576715774866 & 0.00051\\
10 & 3.141592646875561 & 0.00000\\
15 & 3.141592653585951 & 0.00000\\
20 & 3.141592653589791 & 0.00000\\
25 & 3.141592653589793 & 0.00000\\
30 & 3.141592653589793 & 0.00000\\
35 & 3.141592653589793 & 0.00000\\
40 & 3.141592653589793 & 0.00000\\
45 & 3.141592653589793 & 0.00000\\
50 & 3.141592653589793 & 0.00000\\
55 & 3.141592653589793 & 0.00000\\
60 & 3.141592653589793 & 0.00000\\
65 & 3.141592653589793 & 0.00000\\
70 & 3.141592653589793 & 0.00000\\
75 & 3.141592653589793 & 0.00000\\
80 & 3.141592653589793 & 0.00000\\
85 & 3.141592653589793 & 0.00000\\
90 & 3.141592653589793 & 0.00000\\
95 & 3.141592653589793 & 0.00000\\
100 & 3.141592653589793 & 0.00000\\
1000    & 3.141592653589793 & 0.00000\\ 
10000   & 3.141592653589793 & 0.00000 \\ 
100000 & 3.141592653589793 & 0.00000\\ 
1000000  & 3.141592653589793 & 0.00000\\ 
10000000  & 3.141592653589793 & 0.00000\\ [1ex]
\hline
\end{tabular}
\end{table}

\newpage

\subsection{Continued Fractions for $\pi$}
\begin{table}[ht]
\caption{As $n$ increases how does the length of fractions converge  towards $\pi$}
\centering
\begin{tabular}{c c c}
\hline\hline 
n value & Method\#4 & Error(\%) \\ [0.5ex]
\hline
1 & 2.666666666666667 & 15.11741 \\
2 & 2.800000000000000 & 10.87322 \\
3 & 3.428571428571428 & 9.134824 \\
4 & 2.911111111111111 & 7.33646 \\
5 & 3.331601731601732 & 6.04818 \\
6 & 2.980708180708181 & 5.12111 \\
7 & 3.280808080808081 & 4.43136 \\
8 & 3.019032601385542 & 3.90121 \\
9 & 3.250989945726788 & 3.48222 \\
10 & 3.042842125195067 & 2.11028 \\
15 & 3.207889026381334 & 1.58658 \\
20 & 3.091748884841698 & 1.27069\\
25 & 3.181512659824787 & 1.05956\\
30 & 3.108305614026907 & 0.90853\\
35 & 3.170134928816534 & 0.79515\\
40 & 3.116612184235621 & 0.70692\\
45 & 3.163801158726882 & 0.63630\\
50 & 3.121602653391091 & 0.56818\\
55 & 3.123742628422234 & 0.53033\\
60 & 3.124931773880152 & 0.48956\\
65 & 3.156972717366711 & 0.45388\\
70 & 3.127333523307754 & 0.42432\\
75 & 3.154923023986425 & 0.39781\\
80 & 3.129095094977018 & 0.37442\\
85 & 3.153355324069958 & 0.35362 \\
90 & 3.130483257145759 & 0.33502\\
95 & 3.152117511448855 & 0.31827\\
100 & 3.131593903583553 & 0.29872\\ [1ex]
\hline
\end{tabular}
\end{table}

\newpage

\subsection{Viete's Formula of Nested Fractions for $\pi$}
\begin{table}[ht]
\caption{As $n$ increases how does Vietes formula converge  towards $\pi$}
\centering
\begin{tabular}{c c c}
\hline\hline 
n value & Method\#5 & Error(\%) \\ [0.5ex]
\hline
1 & 3.061467458921242 & 2.55046 \\
2 & 3.121445152263491& 0.64135 \\
3 & 3.136548490541725 & 0.16056 \\
4 & 3.140331156952385 & 0.04015\\
5 & 3.141277250932773 & 0.01004 \\
6 & 3.141513801175428 & 0.00251 \\
7 & 3.141572940255612 & 0.00063\\
8 & 3.141587725373528 & 0.00016\\
9 & 3.141591413562714 & 0.00004\\
10 & 3.141592345570118 & 0.00000\\
15 & 3.141592653288993 & 0.00000\\
20 & 3.141592653589493 & 0.00000\\
25 & 3.141592653589793 & 0.00000\\
30 & 3.141592653589793 & 0.00000\\
35 & 3.141592653589793 & 0.00000\\
40 & 3.141592653589793 & 0.00000\\
45 & 3.141592653589793 & 0.00000\\
50 & 3.141592653589793 & 0.00000\\
55 & 3.141592653589793 & 0.00000\\
60 & 3.141592653589793 & 0.00000\\
65 & 3.141592653589793 & 0.00000\\
70 & 3.141592653589793 & 0.00000\\
75 & 3.141592653589793 & 0.00000\\
80 & 3.141592653589793 & 0.00000\\
85 & 3.141592653589793 & 0.00000\\
90 & 3.141592653589793 & 0.00000\\
95 & 3.141592653589793 & 0.00000\\
100 & 3.141592653589793 & 0.00000\\ [1ex]
\hline
\end{tabular}
\end{table}

\newpage

\subsection{The Zeta Function}
\begin{table}[ht]
\caption{As $n$ increases how do the different Zeta Series converge towards $\pi$}
\centering
\begin{tabular}{c c c c c}
\hline\hline 
n value & $\zeta(2)$ & $\zeta(4)$ & $\zeta(6)$ & $\zeta(8)$ \\ [0.5ex]
\hline
5 & 3.09466952411370 & 3.14016117947426 & 3.14157300346359 & 3.14159231269578\\
10 & 3.04936163598207 & 3.14138462246697 & 3.14159185608168 & 3.14159264970117\\
15 & 3.07938982603209 & 3.14152783068467 & 3.14159253913011 & 3.14159265333235\\
20 & 3.09466952411370 & 3.14156460959141 & 3.14159262524305 & 3.14159265355327\\
25 & 3.10392339170058 & 3.14157807684660 & 3.14159264406125 & 3.14159265358185\\
30 & 3.11012872814126 & 3.14158413278489 & 3.14159264969505 & 3.14159265358752\\
35 & 3.11457886229313 & 3.14158724909022 & 3.14159265176594 & 3.14159265358901\\
40 & 3.11792619829938 & 3.14158901347572 & 3.14159265264583 & 3.14159265358948\\
45 & 3.12053546308708 & 3.14159008631043 & 3.14159265306227 & 3.14159265358966\\
50 & 3.12262652293373 & 3.14159077577492 & 3.14159265327654 & 3.14159265358973\\
55 & 3.12433980504914 & 3.14159123889581 & 3.14159265339440 & 3.14159265358976\\
60 & 3.12576920214052 & 3.14159156142938 & 3.14159265346284 & 3.14159265358977\\
65 & 3.12697987310384 & 3.14159179291850 & 3.14159265350444 & 3.14159265358978\\
70 & 3.12801845342065 & 3.14159196334879 & 3.14159265353070 & 3.14159265358979\\
75 & 3.12891920064047 & 3.14159209159432 & 3.14159265354784 & 3.14159265358979\\
80 & 3.12970784547462 & 3.14159218993944 & 3.14159265355935 & 3.14159265358979\\
85 & 3.13040408931831 & 3.14159226661411 & 3.14159265356727 & 3.14159265358979\\
90 & 3.13102327252367 & 3.14159232727297 & 3.14159265357284 & 3.14159265358979\\
95 & 3.13157751780151 & 3.14159237588858 & 3.14159265357684 & 3.14159265358979\\
100 & 3.13207653180911 & 3.14159241530737 & 3.14159265357975 & 3.14159265358979 \\[1ex]
\hline
\end{tabular}
\end{table}

\newpage

\begin{table}[ht]
\caption{As $n$ increases how does the percentage error converge towards zero}
\centering
\begin{tabular}{c c c c c}
\hline\hline 
n value & $\zeta(2)$ & $\zeta(4)$ & $\zeta(6)$ & $\zeta(8)$ \\ [0.5ex]
\hline

5 & 1.49361 & 0.04557 & 0.00041 & 0.00001 \\
10 & 2.93582 & 0.00662 & 0.00017 & 0.00000\\
15 & 1.97998 & 0.00206 & 0.00005 & 0.00000\\
20 & 1.49361 & 0.00089 & 0.00003 & 0.00000\\
25 & 1.19905 & 0.00046 & 0.00000 & 0.00000\\
30 & 1.00153 & 0.00027 & 0.00000 & 0.00000\\
35 & 0.85988 & 0.00017 & 0.00000 & 0.00000\\
40 & 0.75333 & 0.00012 & 0.00000 & 0.00000\\
45 & 0.67027 & 0.00008 & 0.00000 & 0.00000\\
50 & 0.60371 & 0.00006 & 0.00000 & 0.00000\\
55 & 0.54918 & 0.00003 & 0.00000 & 0.00000\\
60 & 0.50368 & 0.00001 & 0.00000 & 0.00000\\
65 & 0.46514 & 0.00000 & 0.00000 & 0.00000\\
70 & 0.43208 & 0.00000 & 0.00000 & 0.00000\\
75 & 0.40341 & 0.00000 & 0.00000 & 0.00000\\
80 & 0.37831 & 0.00000 & 0.00000 & 0.00000\\
85 & 0.35614 & 0.00000 & 0.00000 & 0.00000\\
90 & 0.33643 & 0.00000 & 0.00000 & 0.00000\\
95 & 0.31879 & 0.00000 & 0.00000 & 0.00000\\
100 & 0.30291 & 0.00000 & 0.00000 & 0.00000\\ [1ex]
\hline
\end{tabular}
\end{table}

\newpage
\includepdf[pages={1,2,3,4,5,6,7,8,9}]{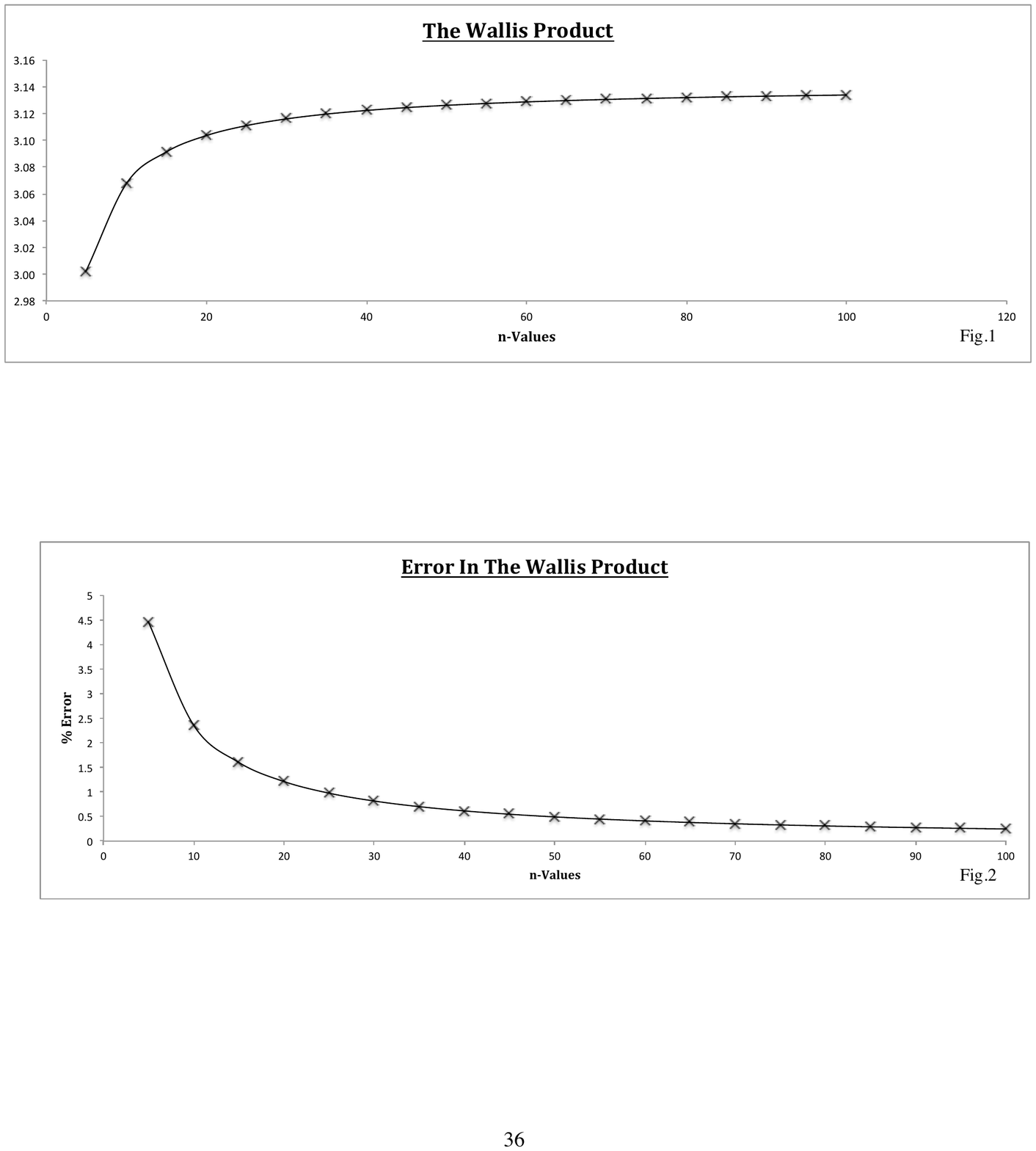}

\section{Conclusions}
\subsection{Conclusion: The Wallis Product}
The Wallis product provides a partially slow method of convergence towards $\pi$. It can be seen from the data table and percentage error that the Wallis Product does not oscillate between two values but instead converges at a slower rate shown by the shallowness of the curve of Fig.1. The Values of $\pi$ generated are always increasing therefore the graph attained can be deemed plausible and correct as it curves upward in Fig.1 . The percentage error shows that even when $n = 10000000$ the actual value of $\pi$ attained contains an error at only the 7th decimal place.
\\ \\
In addition, the ${\%}$ error in Fig.2 gives rise to a curve which starts off with a steep fairly high percentage error and which decreases rapidly and then the rate at which the values attained become more accurate is slowed as the curve becomes much more shallow. On inspection of the error curve from n=5 and n=10 with a comparison to n=55 and n=60 this hypothesis can be proven correct.
\\ \\
An explanation of this lies within the fact that the product only has small increments of numbers for example:

\begin{table}[ht]
\centering
\begin{tabular}{c c}
\hline\hline 
n value & $\frac{2n}{2n-1} \times \frac{2n}{2n+1}$ \\ [0.5ex]
\hline
3 & $\frac{36}{35}$ \\ \\
4 & $\frac{64}{63}$ \\ \\
5 & $\frac{100}{99}$ \\ [1ex]
\hline
\end{tabular}
\end{table}

Here, it can be seen that the fact that the Wallis product converges so slowly is due to the fact that it is being increased  by such minor values that it would require a large amount of iterations to approach the true value of $\pi$ at 15 decimal places

\subsection{Conclusion: Eulers Continued Fractions}

The Continued fraction used in this method for calculating $\pi$ is one of the slowest continued fractions formula used to calculate pi. An explanation for this inefficiency (Fig.9) lies in the fact that the values of $\pi$ per iteration alternate and oscillate between the positive and negative regions about the line $y = \pi$ very slowly and in turn this results for a very slow convergence.
\\ \\ 
By inspecting the data from the table and the graph in Fig.10, the percentage error's clearly illustrate how inefficient and inappropriate the method is for calculating accurate values of $\pi$ as when n = 1 the percentage error is quite significant with a large 15.1{\%}. Furthermore, the gradient and drop in percentage error is not as fast as the other methods which accounts for the very shallow curve of convergence towards 0{\%}.
\\ \\
If we look at the actual continued fraction formula itself,  the value of each section grows very slowly and therefore has a very slow impact on the value of the whole continued fraction itself. This accounts for part of the inefficiency of this method. In addition, the continued fraction that Euler developed his so slow that it requires ''roughly $3 \times 10^{n}$ terms to achieve n-decimal precision" \cite{cite_key11}.

\subsection{Conclusion: The Gregory-Leibniz Series}
The Gregory-Leibniz Series illustrates an out of date and inefficient infinite series which will converge to $\pi$ at a very slow rate. Some partial reasoning behind this is the fact that the actual infinite sum contains $(-1)^n$in the numerator and therefore causes the partial sums to oscillate between the positive and negative regions above and below the line $y=\pi$. This can be seen through Fig.3 that the attained $\pi$ values which alternate between every $n$ value either going far much above the required value or too far below.
\\ \\
Both the Table of Data and the graph emphasize the inefficiency of the infinite sum. Fig.3 illustrates the slow convergence towards $\pi$ as the graph is not very steep at the start and becomes even more shallow as the n values increase which is also supported by the fact that "Calculating $\pi$ to 10 correct decimal places using direct summation of the series requires about 5,000,000,000 terms"\cite{cite_key4} .
\\ \\
In inspecting the inefficiency, it is important to note that that even though Mathematica would display a percentage error of 0.00000{\%} at the $10000000^{th}$ $n$ term there is still a margin of error as the infinite series makes an error at only the 8th decimal point. Taking into account that it was discovered in 1670's where mathematical analysis was not yet fully developed and that it was the "first ever found infinite series for $\pi$"\cite{cite_key4} , It would have been a tremendous discovery by James Gregory and Gottfried Wilhelm Von Leibniz at the time. Furthermore, "In order to achieve 100 accurate decimal places, one would have to go through $1\times10^{47}$ iterations of the series" \cite{cite_key4} which illustrates how it cannot satisfy modern needs for the computation and calculation of accurate values of $\pi$ shown in Fig.4.

\subsection{Conclusion: Newtons Arcsine Series Expansion}

Newton's Arcsine Expansion provides a rapidly converging infinite series which attains an accurate value of $\pi$ between n=20 and n=10. An explanation for this accelerated convergence can be seen within the actual infinite function. If one pays close attention to the components of the infinite series, it can be seen that $(2n)!$ and $(n!)^2$ are present on the numerator and denominator respectively. The actual values of these factorials expand rapidly as $n$ increases. This can be shown in the sub-table below:
\\ \\

\begin{table}[ht]
\centering
\begin{tabular}{c c c}
\hline\hline 
n value & $(2n)!$ & $(n!)^2$ \\ [0.5ex]
\hline
3 & 720 & 36\\
4 & 40320 & 576\\
5 & 3628800 & 14400 \\ [1ex]
\hline
\end{tabular}
\end{table}

Of course the value of (2n)! will increase at a much more rapid rate which also converges to a limit fold which would further then explain the shape of the convergence graph as curving upwards exponentially as it  approaches  the value of $\pi$ in Fig.5. Also, not that the middle binomial coefficient is present in Newton's Arcsine expansion which is $\begin{pmatrix} 2n \\ n \end{pmatrix}$ and this also has a high numerical value as n becomes large which accounts for the rapid convergence.
\\ \\
In addition, by inspecting the both the data tables and Fig.6 for Newton's Arcsine expansion, it can be seen that between the terms in which n=5 and n=10 the series assumes a percentage error of less that 0.00005{\%} which accounts for the extremly steep slope in the graph for percentage error and the extremely fast convergence in the graph displaying the corresponding $\pi$ values. Alternatively, Sterling's approximation could also be used here to investigate the rate of convergence.

\subsection{Conclusion: Vieta's Nested Radical Formula}
Vieta’s Formula consists of nested radicals which provide a fairly fast convergence towards $\pi$. From the table, it can be seen that at the start with the first iteration at the first n value there is approximately 2.55{\%} error and provides an incorrect value at the first decimal point. However it can be seen through the steepness of the percentage error in Fig.8 between $n=1$ and $n=2$ that the error rapidly decreases which can be seen by thorough inspection at the second percentage error where n=2 has decreased by approximately a factor of 4. Furthermore, by considering the n value in which Vieta's Formula provides an accurate value of $\pi$ it can be seen that between $n=10$ and $n=15$ the re-iterative values reach the correct display of $\pi$ for 15 decimal places which remains consistent for the rest of the iterations.
\\ \\
Furthermore, the explanation behind this lies within the function itself. The whole use of nested radicals allows the function to not fluctuate between the positive and negative regions about the line $y =\pi$ and only use the positive region. Furthermore, if we consider the following part of Vieta's Formula:

\begin{eqnarray}
\sqrt{2-\underbrace{{\sqrt{2+\sqrt{2+\sqrt2}}}}_a} = a \nonumber
\end{eqnarray}

Then, by squaring both sides of the equation,
\begin{eqnarray}
2 + a &=& a^2 \nonumber \\
a^2 - a - 2 &=& 0 \nonumber \\
(a-2)(a+1) &=& 0 \nonumber
\end{eqnarray}

Since $a$ must be positive, $a = 2$ 

From this, if we now consider the whole form of Vieta's Formula:

\begin{eqnarray}
2^{2n+1} \sqrt{2-\underbrace{{\sqrt{2+\sqrt{2+\sqrt2}}}}_{a \rightarrow 2}...}\nonumber
\end{eqnarray}

Therefore, as $n$ increases the value of $2^{2n+1}$ will increase and the value of $a$ nested roots two radicals will approach two and the value of the square root will decrease. However, it is important to note that the speed at which $2^{2n+1}$ increases will be faster than the speed that the nested roots approach 2. This in turn explains the shape of the graph as it provides evidence that the graph should rapidly accelerate upwards because at first the rate at which $2^{2n+1}$ increases will be much greater than the rate of decrease in the nested radicals. Then, as shown in fig.7, the curve should flatten as the rates increase of both terms equal out until they are balanced at which n tends towards infinity ${\infty}$ and the exact value of ${\pi}$ is obtained.

\subsection{Conclusion: The Reimann-Zeta Function}
The series developed by the Reimann-Zeta Function for $\zeta(2),\zeta(4),\zeta(6),\zeta(8)$ show an increase in accuracy as the the value of $n$ in the Zeta function increases. This in turn provides evidence for the attained graphs for the different $n$ values. For example, by looking at Fig. 11 for $\zeta(2)$ one can notice that the gradient flattens out at a much slower rate than which accounts for a slower rate of convergence compared to Fig.13 where there is a much steeper gradient at the start and the rate at which the series converges is much faster.
\\ \\
Furthermore, by considering the various graphical displays for the error values, one can notice as expected that as the value of $n$ increases, the amount of error in the series decreases at a much faster rate. Take for example the error in Fig.15 with the graph of $\zeta(6)$, Here it is evident that the gradient at the start is quite steep and begins to touch close to 0\% at about $n=20$ Similarly, if one examine's the graph of error in $\zeta(8)$ through Fig.18 it can be seen that the gradient is much steeper at the start which signifies how accurate this method is in calculating $\pi$.
\\ \\
The explanation behind this lies in the Zeta function itself. Firstly, consider how the nature of how $\zeta(2)$ is defined. The series is based on having a square root of the summation of series with fractions with a numerator of 1 and denominator of $n^2$. Because of this, each additional fraction added in the series will have less of an influence on the value of $\pi$ attained and in turn part of the series accounts for a slower convergence towards $\pi$ in Fig.11. Contrastingly, if one pays close attention to the series obtained by $\zeta(8)$, it is evident that here, the square root is raised to a power of 8 and each additional fraction added in the summation has a numerator of 1 and denominator of $n^8$. Because of this, each additional fraction added becomes even less influential to the value of $\pi$ attained which can be seen in the table below. Conclusively, this accounts for the shape of the error curves in the Figures 12 and 18 as Fig.12 is takes an $n$ value above $n=100$ to drop below 0.3\% error but for $\zeta(8)$, between $n=5$ and $n=10$ the\% error already drops to 0.0\% which shows a much steeper curve in Fig.18.

\begin{table}[ht]
\centering
\begin{tabular}{c c}
\hline\hline 
$\zeta(n)$ & 3rd Fractional Order Value \\ [0.5ex]
\hline
$\zeta(2)$ & $\frac{1}{3^2}=0.11111$\\
$\zeta(4)$ & $\frac{1}{3^4}=0.01232$\\
$\zeta(8)$ & $\frac{1}{3^6}=0.00137$\\
$\zeta(6)$ & $\frac{1}{3^8}=0.00002$\\ [1ex]
\hline
\end{tabular}
\end{table}

\section{Comparisons of the different methods}
\subsection{The Gregory-Leibniz Series VS Newtons Series: The Arcsine Expansion}
\begin{eqnarray}
\sum_{n=0}^{\infty} \frac{(-1)^n}{2n+1} x^{2n+1}    \mbox{			And			}        \sum_{n=0}^{\infty} \frac{(2n)!}{(2n+1)(n!)2^{2n}} x^{2n+1} \nonumber
\end{eqnarray}

The infinite series on the left represents the Gregory-Leibniz series and the infinite series on the right represents Newtons infinite series for the arcsine function. By paying close attention to the terms in common and only considering the terms that are not in common we obtain the following (as $x=1$ in the Gregory-Leibniz Series and $x=0.5$ in Newton's Arcsine Expansion:

\begin{eqnarray}
\sum_{n=0}^{\infty} (-1)^n  \mbox{			And			}        \sum_{n=0}^{\infty} \frac{(2n)!}{(n!)2^{2n}} \nonumber
\end{eqnarray}

Considering this, it is obvious that the second infinite sum will converge much faster than the first because as previously stated in the analysis for Newtons infinite series, the value of $\frac{(2n)!}{2^{2n} (n!)}$ will converge at a much faster rate because the use of $(2n!)$ and $(n!)$ allows the series to converge towards a value much faster than the infinite series containing $(-1)^n$ because there is no fluctuation.
\\ \\
Furthermore, by having a close look at Fig.19, The graph of comparison clearly confirms the fact that Newton's arcsine expansion from previous reasoning is a much more efficient way to achieve an accurate value of pi as by n=30, the Newton's arcsine expansion curve is practically linear whilst the Gregory-Leibniz oscillations are still oscillating in wide amounts.
\\ \\
Also, Fig.20 yet again emphasizes how accurate and efficient the arcsine expansion is compared to the Gregory-Leibniz series as the red plot is always underneath the blue plot. It is important to note that there is some percentage error and the line $y=0$ is an asymptote as there will always be some error however Mathematica was not able to calculate to that precision.

\subsection{Vieta's Nested Radicals VS Eulers Continued Fraction}

\begin{eqnarray}
\lim_{n \rightarrow \infty} 2^{n+1} \sqrt{\sqrt 2- \underbrace{\sqrt {2+\sqrt{2 + \sqrt{2 + ... + \sqrt2}}}}_{n}} \mbox{ 			And			} \cfrac{4}{1^2+\cfrac{1^2}{2+\cfrac{3^2}{2+\cfrac{5^2}{2+\cfrac{7^2}{2+\cfrac{9^2}{\ddots}}}}}} \nonumber 
\end{eqnarray}

Vieta's Nested radicals provide a much faster convergence than Eulers Continued fraction for $\pi$ because Euler's continued fraction provides a very slow increase in numerators of the continued fractions as they increase by the unique sequence of odd numbered squares. By considering how the nested $\sqrt2$ radicals work we can see that the fact that as n = 3 the nested fractions are already near complete convergence at 1.96157 and approaching towards 2 which proves why nested radicals are more efficient in calculating $\pi$
\\ \\
Using graphical evidence, Fig.22 provides evidence for the previously mentioned hypothesis as said before Vieta's nested radicals converges at a much faster rate than Eulers continued Fractions. While Vieta's continued fractions have reached n=25, Euler's continued fraction is still oscilating above and below the line $y=\pi$ Also, Fig.21 illustrates how there is virtually no percent error for Vieta's nested fractions past n=25 but at n=25 Euler's continued fractions still obtains a percentage error of 1.5${\%}$. Also, it becomes more difficult to calculate the actual values for Vieta's square roots because each additional square root increase the likelyhood for a rounding error.

\subsection{The Wallis Product VS Newtons Series: The Arcsine Expansion}

\begin{eqnarray}
\prod^{\infty}_{k=1} 2\left[\frac{2k}{2k-1} . \frac{2k}{2k+1}\right] \mbox{ 			And 				} \sum_{n=0}^{\infty} \frac{(2n)!}{(2n+1)(n!)2^{2n}} x^{2n+1} \nonumber 
\end{eqnarray}
 
The product on the left represents the wallis product for pi and the infinite series on the right represents Newtons Arcsine expansion. Here it is clear that Newtons Arcsine expansion will give a much faster convergence rate. Take the transition between $k=6$ $\rightarrow$ $k=7$ and $n=6$ to $n=7$ as examples. The value of the numerator of the wallis product will increase by 104 (See Calculation 1.00) and the denominator by a factor of 13x15 which is a minor increment compared to the numerator in the Arcsine function which increases by $8.67 \times 10^{10}$. The power of the factorial allows for a faster convergence then by multiplication hence the reason why the infinite series will converge faster than the product.

\begin{eqnarray}
[1.00] \mbox{	change} = 2.\frac{2.2.4.4.6.6.8.8.10.10.12.12}{1.3.3.5.5.7.7.9.9.11.11.13}.(\frac{14.14.}{13.15}-1) \nonumber
\end{eqnarray}

As shown in Fig.23, The speed of convergence for Newton's Arcsine still remains unparalleled as it is able to converge much more rapidly than the Wallis product as the blue plot is linear and the purple plot remains as a curve. Fig.24 also illustrates the accuracy of Newton's arcsine expansion as the percentage error values are constantly far below the percentage error values of the Wallis Product

\subsection{The Wallis Product VS The $\zeta(2)$ Series}
\begin{eqnarray}
\prod^{\infty}_{k=1} 2\left[\frac{2k}{2k-1} . \frac{2k}{2k+1}\right] \mbox{ 			And 				} \sum^{\infty}_{n=0} \sqrt{6\left[\frac{1}{1^2}+\frac{1}{2^2}+\frac{1}{3^2}+\frac{1}{4^2}+...\right]} \nonumber
\end{eqnarray}

The product on the left is the Wallis Product for $\pi$ and the infinite series on the right represents the infinite series created at $\zeta(2)$ which in turn converges towards $\pi$. Here it is clear that the series for $\zeta(2)$ converges at a slightly slower rate than the Wallis Product because of the fact that firstly the infinite series makes use of multiplication of 2 instead whilst the Wallis Product have a larger multiplication increment between each n value which is greater than 2 and secondly because of the fact that a square root is used in place of multiplication which causes the $\zeta(2)$ series to be slower.
\\ \\

By examining the graphical display in Fig.25, it is clear that the Wallis Product will converge at a faster rate which is shown as the graph representing the $\zeta(2)$ infinite series is always slightly below the graph of the Wallis Product. Furthermore, by taking into consideration the $\%$ error in each of the methods, the results shown in Fig.26 consolidate what was previously mentioned as at all values for $n$, the graph of the infinite series with $\zeta(2)$ is always above the graph of the Wallis Product which in turn shows that the rate of accuracy increases at a much larger rate in the Wallis product. Also, by taking into account the gradient of the percentage error curves, between values $n=5$ and $n=10$ one can see that the gradient of the line for the Wallis Product is much steeper than the gradient of the line for $\zeta(2)$ which emphasizes the fact that the rate of accuracy for the Wallis Product increases at a slightly faster pace than the rate of accuracy for the infinite series.

\subsection{Newton's Arcsine Expansion VS The  $\zeta(8)$ Series}
\begin{eqnarray}
\sum_{n=0}^{\infty} \frac{(2n)!}{(2n+1)(n!)2^{2n}} x^{2n+1} \nonumber  \mbox{                And                  }
\sum^{\infty}_{n=0} \sqrt[8]{9450\left[\frac{1}{1^8} + \frac{1}{2^8} +\frac{1}{3^8} + \frac{1}{4^8} + ... \right] \nonumber}
\end{eqnarray}

The infinite series on the left represents one of the most efficient ways to attain a value of $\pi$ developed by Newton and his expansion of arcsine. The series on the right is the fastest way in attaining a value of $\pi$ covered in this thesis which is the infinite series formed using $\zeta(8)$. In this case, it is clear that not only does  Newton's arcsine expansion makes use of rapidly increasing factorial values shown in $(2n)!$ and $(n!)$ but also uses high valued exponentials such as $2^{2n}$. However, these factor's are outweighed by the infinite series developed from the Reimann-Zeta Function of value $\zeta(8)$ because here there is not only multiplication by 9450 but also makes use of the 8th square root and every proceeding term in the summation begins to affect the final product less and less at a much faster rate that the summation of Newton's Arcsine Expansion.
\\ \\
Furthermore, if we consider the graphical display in Fig.27 and Fig.28 at $n=5$, the produced value the $\zeta(8)$ series where $n=5$ is much closer to $\pi$ than the produced value produced by Newtons Arcsine Expansion. Also, The \% error where $n=5$ for $\zeta(8)$ is 26 times smaller than the produced value than the \% error given through newtons arcsine expansion which emphasizes the degree of accuracy that the infinite series that $\zeta(8)$ has over the series from Newton's Arcsine expansion.

\newpage

\includepdf[pages={1,2,3,4}]{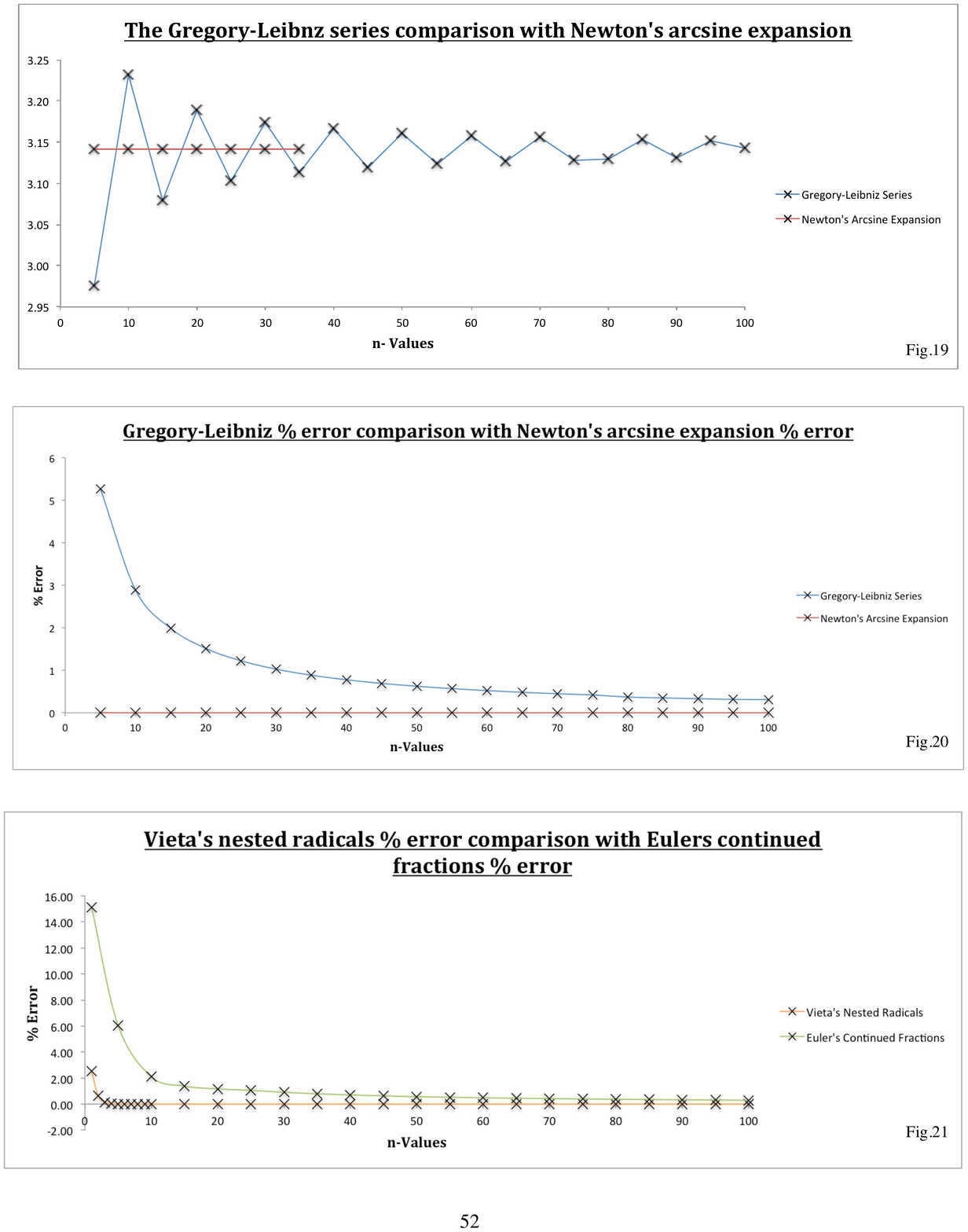}

\section{Appendix}
\subsection{Mathematic Iteration code: The Wallis Product}

"Wallis Product"
For[i = 0, i <= 100, i = i + 5, Print[i " SymbA ", NumberForm[N[2 /!/(
/*UnderoverscriptBox[/(/[Product]/), /(k = 1/), /(i/)]/((
/*FractionBox[/(2  k/), /(2  k - 1/)] 
/*FractionBox[/(2  k/), /(2  k + 1/)])/)/)], 16]]]

5  SymbA 3.002175954556907

10 SymbA 3.067703806643499

15  SymbA 3.091336888596228

20  SymbA 3.103516961539234

25  SymbA 3.11094516690154

30  SymbA 3.115948285887959

35  SymbA 3.119547206305518

40  SymbA 3.122260326421437

45  SymbA 3.124378835915516

50  SymbA 3.126078900215411

55  SymbA 3.127473350412857

60  SymbA 3.128637797891591

65  SymbA 3.129624812079802

70  SymbA 3.130472076319065

75  SymbA 3.131207308587379

80  SymbA 3.131851351372613

85  SymbA 3.132420179022906

90  SymbA 3.132926240627509

95  SymbA 3.133379381619937

100 SymbA 3.133787490628162

1000 SymbA 3.140807746030395

10000 SymbA 3.141514118681922

100000 SymbA 3.141584799657247

1000000 SymbA 3.141591868192124

10000000 SymbA 3.141592575049982
\subsubsection{Mathematic Iteration code: The Wallis Product Percentage Error}

\begin{multicols}{2}
\begin{itemize}
    \item 1 - N[{3.002175954556907/π}, 15]
N[{0.044377713601277735`*100}, 15]
4.43777
    \item N[{(3.067703806643499/π)}, 15]
1 - {0.9764804495382736`}
{0.02351955046172638`}
*100 = 2.35196
    \item 3.091336888596228/π
1 - 0.9840030931648189`
0.01599690683518107`*100
1.59969

    \item 3.103516961539234/π
1 - 0.9878801307970174`
0.01211986920298258*100
1.21199

    \item 3.110945166901554/π
1 - 0.9902446020004474`
100*0.009755397999552606`
0.97554

    \item 3.115948285887959/π
1 - 0.9918371442355739`
100*0.008162855764426102`
0.81629

\item 1 - 0.9929827161840716`
100*0.007017283815928366`
0.701728

\item 3.122260326421437/π
1 - 0.9938463291393728`
100*0.006153670860627236`
0.615367

\item 3.124378835915516/π
1 - 0.9945206716553123`
100*0.005479328344687673`
0.547933

\item 3.126078900215411/π
1 - 0.9950618189291169`
100*0.004938181070883063`
0.493818

\item3.127473350412857/π
1 - 0.9955056862127551`
100*0.004494313787244852`
0.449431

\item1 - 3.128637797891591/π
100*0.004123658
642822159` = 0.412366

\item3.129624812079802/ π
1 - 0.996190517731089`
100*0.0038094822689109797`
0.380948

\item3.130472076319065/π
1 - 0.9964602103146565`
100*0.0035397896853435196`
0.353979

\item3.131207308587379/π
1 - 0.9966942420143022`
100*0.003305757985697766`
0.330576

\item3.131851351372613/π
1 - 0.9968992471999675`
100*0.0031007528000325246`
0.310075

\item3.132420179022906/π
1 - 0.9970803106645905`
100*0.0029196893354095366`
0.291969

\item3.132926240627509/π
1 - 0.9972413950763536`
100*0.002758604923646435`
0.27586

\item 3.133379381619937/π
1 - 0.9973856343340787`
100*0.0026143656659213255`
0.261437

\item 3.133787490628162/π
1 - 0.9975155394660373`
100*0.0024844605339626735`
0.248446

\item 3.140807746030395/π
1 - 0.9997501561641032`
100*0.00024984383589676806`
0.0249844

\item 3.141514118681922/π
1 - 0.9999750015624141`
100*0.000024998437585854738`
0.00249984

\item 3.141584799657247/π
1 - 0.9999975000156252`
100*2.4999843748485517`*
0.000249998

\end{itemize}
\end{multicols}

\subsection{Mathematic Iteration code: The Gregory-Leibniz Series}

For[i = 0, i <= 100, i = i + 5, Print[i "  SymbA  ", NumberForm[N[4 /!(/
/*UnderoverscriptBox[/(/[Sum]/), /(n = 0/), /(i/)]
/*FractionBox[
SuperscriptBox[/((/(-1/))/), /(n/)], /(2  n + 1/)]/)], 16]]]

5  SymbA 2.976046176046176

10 SymbA 3.232315809405593

15 SymbA 3.079153394197426

20 SymbA  3.189184782277595

25 SymbA  3.103145312886011

30 SymbA  3.173842337190749

35 SymbA  3.113820229023573

40 SymbA  3.165979272843215

45 SymbA  3.119856090062712

50 SymbA  3.161198612987056

55 SymbA  3.123736933726277

60 SymbA  3.157984995168666

65 SymbA  3.126442007766234

70 SymbA  3.155676462307475

75 SymbA  3.128435328236984

80 SymbA  3.153937862272616

85 SymbA  3.129965139593801

90 SymbA  3.152581332875124

95 SymbA  3.131176269454981

100 SymbA  3.151493401070910

1000 SymbA 3.142591654339543

10000 SymbA 3.141692643590543

100000 SymbA 3.141602653489794

1000000 SymbA 3.141591868192127

10000000 SymbA  3.141592653518272

\newpage

\subsubsection{Mathematic Iteration code: The Gregory-Leibniz Series Percentage Error}

\begin{multicols}{2}
\begin{itemize}
\item 2.976046176046176/π
1 - 0.9473049195749638`
100*0.052695080425036234`
5.26951

\item 3.232315809405593/π
1 - 1.0288780774019617`
100*-0.028878077401961688`
-2.88781

\item 3.079153394197426/π
1 - 0.9801249664494155`
100*0.01987503355058451`
1.9875

\item 3.189184782277595/π
1 - 1.0151490450658587`
100*-0.01514904506585868`
-1.5149

\item 3.103145312886011/π
1 - 0.9877618313565096`
100*0.012238168643490366`
1.22382

\item 3.173842337190749/π
1 - 1.0102653931164836`
100*-0.010265393116483557`
-1.02654

\item 3.113820229023573/π
1 - 0.9911597626972786`
100*0.008840237302721432`
0.884024

\item 3.165979272843215/π
1 - 1.0077625019989642`
100*-0.007762501998964222`
-0.77625

\item 3.119856090062712/π
1 - 0.993081036937668`
100*0.006918963062331973`
0.691896

\item 3.161198612987056/π
1 - 1.0062407707042669`
100*-0.006240770704266874`
-0.624077

\item 3.123736933726277/π
1 - 0.9943163478425145`
100*0.005683652157485453`
0.568365

\item 3.157984995168666/π
1 - 1.0052178443822568`
100*-0.005217844382256809`
-0.521784

\item 3.126442007766234/π
1 - 0.9951773996522919`
100*0.004822600347708095`
0.48226

\item 3.155676462307475/π
1 - 1.0044830155499596`
100*-0.004483015549959557`
-0.448302

\item 3.128435328236984/π
1 - 0.9958118932644642`
100*0.004188106735535824`
0.418811

\item 3.129965139593801/π
1 - 0.9962988473433353`
100*0.003701152656664708`
0.370115

\item 3.152581332875124/π
1 - 1.0034978052526238`
100*-0.0034978052526237757`
-0.349781

\item 3.131176269454981/π
1 - 0.9966843619516014`
100*0.003315638048398628`
0.331564
\item 3.151493401070910/π
1 - 1.003151505803849`
100*-0.0031515058038489308`
-0.315151

\item 3.142591654339543/π
1 - 1.0003179918149503`
100*-0.0003179918149502914`
-0.0317992

\item 3.141692643590543/π
1 - 1.0000318278057583`
100*-0.000031827805758277705`
-0.00318278

\item 100*-1.00000318306703`*15
1.00000318306703`*100
0.00001

\end{itemize}
\end{multicols}

\subsection{Mathematic Iteration code: Newtons Arcsine Expansion Series}
For[i = 0, i <= 100, i = i + 5, Print[i " SymbA", NumberForm[N[6 /!/(
/*UnderoverscriptBox[/([/Sum]/), /(n = 0/), /(i/)]
/*FractionBox[/(/(/((2  n)/)!/) 
/*SuperscriptBox[/((1/2)/), /(2  n + 1/)]/), /(
/*SuperscriptBox[/(2/), /(2  n/)] 
/*SuperscriptBox[/(((/n!)/)/), /(2/)] /((2  n + 1)/))]/)], 16]]]

5  SymbA 3.141576715774866

10  SymbA 3.141592646875561

15  SymbA 3.141592653585951

20  SymbA 3.141592653589791

25  SymbA 3.141592653589793

30  SymbA 3.141592653589794

35  SymbA 3.141592653589793

40  SymbA 3.141592653589793

45  SymbA 3.141592653589793

50  SymbA 3.141592653589793

55  SymbA 3.141592653589793

60  SymbA 3.141592653589794

65  SymbA 3.141592653589793

70  SymbA 3.141592653589793

75  SymbA 3.141592653589794

80  SymbA 3.141592653589794

85  SymbA 3.141592653589793

90  SymbA 3.141592653589794

95  SymbA 3.141592653589793

100 SymbA 3.141592653589794

1000 SymbA 3.142591654339543

10000 SymbA 3.141692643590543

100000 SymbA 3.141692643590543

1000000 SymbA 3.141692643590543

\subsubsection{Mathematic Iteration code:Newtons Arcsine Series Percentage Error}

\begin{multicols}{2}
\begin{itemize}
\item 3.141592653585951/π
1 - 0.999999999998777`
100*1.2230216839270724`*-12
1.22302*10*-10

\item 3.141592646875561/π
1 - 0.9999999978627935`
100*2.1372065228675297`*-9
2.13721*10*-7

\item 3.141576715774866/π
1 - 0.9999949268359446`
100*5.073164055402479`*-6
0.000507316
\end{itemize}
\end{multicols}

\section{Mathematic Iteration code: Continued Fractions for $\pi$}
g[n,x] = (2 n + 1)*2/(2 + x)
NumberForm[4./(1 + g[0, g[1, g[2, g[3, g[4]]]]]), 16]
‘™
TagBox[
InterpretationBox[""3.331601731601732'",
3.3316017316017317`,
AutoDelete->True],
NumberForm[, 16] ]
NumberForm[4./(
 1 + g[0, g[1, g[2, g[3, g[4, g[5, g[6, g[7, g[8, g[9]]]]]]]]]]), 16]
‘™

TagBox[
InterpretationBox[""3.042842125195067,
3.0428421251950666`,
AutoDelete->True],
NumberForm[, 16] ]
NumberForm[
 4./(1 + g[0, g[1, g[2, g[3, g[4, g[5, g[6, 
           g[7, g[8, g[9, g[10, g[11, g[12, g[13, g[14]]]]]]]]]]]]]]]), 16]
‘™

TagBox[
InterpretationBox[""3.207889026381334",
3.207889026381334,
AutoDelete->True],
NumberForm[, 16] ]
NumberForm[
 4./(1 + g[0, g[1, g[2, g[3, g[4, g[5, g[6,  g[7, g[8,  g[9, g[10,  g[11, g[12, g[13, g[14, g[15, g[16, g[17, g[18, g[19]]]]]]]]]]]]]]]]]]]]), 16]
‘™

TagBox[
InterpretationBox[""3.091748884841698>",
3.0917488848416976`,
AutoDelete->True],
NumberForm[, 16] ]
NumberForm[
 4./(1 + g[0, g[1, g[2, g[3, g[4, g[5, g[6, g[7, g[8, g[9, g[10, g[11, g[12, g[13, g[14, g[15, g[16, g[17, g[18, g[19, g[20, g[21, g[22, g[23, g[24]]]]]]]]]]]]]]]]]]]]]]]]]), 16]
‘™

TagBox[
InterpretationBox[""3.181512659824787",
3.1815126598247874`,
AutoDelete->True],
NumberForm[, 16] ]
NumberForm[
 4./(1 + g[0, 
     g[1, g[2, g[3, g[4, g[5, g[6,  g[7, g[8, g[9, g[10, g[11, g[12, g[13, g[14, g[15, g[16, g[17, g[18, g[19, g[20, g[21, g[22, g[23, g[24, g[25, g[26, g[27, g[28, g[29]]]]]]]]]]]]]]]]]]]]]]]]]]]]]]), 16]
‘™

TagBox[
InterpretationBox[""3.108305614026907",
3.108305614026907,
AutoDelete->True],
NumberForm[, 16] ]
NumberForm[
 4./(1 + g[0, g[1, g[2, g[3, g[4, g[5, g[6, g[7, g[8, g[9, g[10, g[11, g[12, g[13, g[14, g[15, g[16, g[17, g[18, g[19, g[20, g[21, g[22, g[23, g[24, g[25, g[26, g[27, g[28, g[29, g[30, g[31, g[32, g[33, g[34]]]]]]]]]]]]]]]]]]]]]]]]]]]]]]]]]]]), 16]
‘™

TagBox[
InterpretationBox[""3.170134928816534>",
3.1701349288165344`,
AutoDelete->True],
NumberForm[, 16] ]
NumberForm[
 4./(1 + g[0, 
g[1, 
g[2, 
g[3,g[4, g[5, g[6, g[7, g[8, g[9, g[10, g[11, g[12, g[13, g[14, g[15, g[16, g[17, g[18, g[19, g[20, g[21, g[22, g[23, g[24, g[25, g[26, g[27, g[28, g[29, g[30, g[31, g[32, g[33, g[34, g[35, g[36, g[37, g[38, g[39]]]]]]]]]]]]]]]]]]]]]]]]]]]]]]]]]]]]]]]]), 16]
‘™

TagBox[
InterpretationBox[""3.11661218423562",
3.1166121842356196`,
AutoDelete->True],
NumberForm[, 16] ]
NumberForm[
 4./(1 + g[0, g[1, g[2, g[3, g[4, g[5, g[6, g[7, g[8, g[9, g[10, g[11, g[12, g[13, g[14, g[15, g[16, g[17, g[18, g[19, g[20, g[21, g[22, g[23, g[24, g[25, g[26, g[27, g[28, g[29, g[30, g[31, g[32, g[33, g[34, g[35, g[36, g[37, g[38, g[39, g[40, g[41, g[42, g[43, g[44]]]]]]]]]]]]]]]]]]]]]]]]]]]]]]]]]]]]]]]]]]]]]), 16]
‘™

TagBox[
InterpretationBox[""3.163801158726882>",
3.1638011587268817`,
AutoDelete->True],
NumberForm[, 16] ]
NumberForm[
 4./(1 + g[0, g[1, g[2, g[3, g[4, g[5, g[6, g[7, g[8, g[9, g[10, g[11, g[12, g[13, g[14, g[15, g[16, g[17, g[18, g[19, g[20, g[21, g[22, g[23, g[24, g[25, g[26, g[27, g[28, g[29, 
g[30, g[31, g[32, g[33, g[34, g[35, g[36, g[37, g[38, g[39, g[40, g[41, g[42, g[43, g[44, g[45, g[46, g[47, g[48,    g[49]]]]]]]]]]]]]]]]]]]]]]]]]]]]]]]]]]]]]]]]]]]]]]]]]]), 16]
‘™

TagBox[
InterpretationBox[""3.121602653391091>",
3.1216026533910908`,
AutoDelete->True],
NumberForm[, 16] ]
NumberForm[
 4./(1 + g[0, g[1, g[2, g[3, g[4, g[5, g[6, g[7, g[8, g[9, g[10, g[11, g[12, g[13, g[14, g[15, g[16, g[17, g[18, g[19, g[20, g[21, g[22, g[23, g[24, g[25, g[26, g[27, g[28, g[29, g[30, g[31, g[32, g[33, g[34, g[35, g[36, g[37, g[38, g[39, g[40, g[41, g[42, g[43, g[44, g[45, g[46, g[47, g[48, g[49, g[50, g[51, g[52, g[53, g[54, 
g[55]]]]]]]]]]]]]]]]]]]]]]]]]]]]]]]]]]]]]]]]]]]]]]]]]]]]]]]]
), 16]
‘™

TagBox[
InterpretationBox[""3.12374262842223>",
3.1237426284222303`,
AutoDelete->True],
NumberForm[, 16] ]
NumberForm[
 4./(1 + g[0, g[1, g[2, g[3, g[4, g[5, g[6, g[7, g[8, g[9, g[10, g[11, g[12, g[13, g[14, g[15, g[16, g[17, g[18, g[19, g[20, g[21, g[22, g[23, g[24, g[25, g[26, g[27, g[28, g[29, g[30, g[31, g[32, g[33, g[34, g[35, g[36, g[37, g[38, g[39, g[40, g[41, g[42, g[43, g[44, g[45, g[46, g[47, g[48, g[49, g[50, g[51, g[52, g[53, g[54, g[55, g[56, g[57, g[58, g[59]]]]]]]]]]]]]]]]]]]]]]]]]]]]]]]]]]]]]]]]]]]]]]]]]]]]]]]]
]]]]), 16]
‘™

TagBox[
InterpretationBox[""3.12493177388015",
3.1249317738801503`,
AutoDelete->True],
NumberForm[, 16] ]
NumberForm[
 4./(1 + g[0, g[1, g[2, g[3, g[4, g[5, g[6, g[7, g[8, g[9, g[10, g[11, g[12, g[13, g[14, g[15, g[16, g[17, g[18, g[19, g[20, g[21, g[22, g[23, g[24, g[25, g[26, g[27, g[28, g[29, g[30, g[31, g[32, g[33, g[34, g[35, g[36, g[37, g[38, g[39, g[40, g[41, g[42, g[43, g[44, g[45, g[46, g[47, g[48, g[49, g[50, g[51, g[52, g[53, g[54, g[55, g[56, g[57, g[58, g[59, g[60, g[61, g[62, g[63, g[64]]]]]]]]]]]]]]]]]]]]]]]]]]]]]]]]]]]]]]]]]]]]]]]]]]]]]]]]]]]]]]]]]), 16]
‘™

TagBox[
InterpretationBox[""3.156972717366711>",
3.1569727173667115`,
AutoDelete->True],
NumberForm[, 16] ]
NumberForm[
 4./(1 + g[0, g[1, g[2, g[3, g[4, g[5, g[6, g[7, g[8, g[9, g[10, g[11, g[12, g[13, g[14, g[15,g[16, g[17, g[18, g[19, g[20, g[21, g[22, g[23, g[24, g[25, g[26, g[27, g[28, g[29, g[30, g[31, g[32, g[33, g[34, g[35, g[36, g[37, g[38, g[39, g[40, g[41, g[42, g[43, g[44, g[45, g[46, g[47, g[48, g[49, g[50, g[51, g[52, g[53, g[54, g[55, g[56, g[57, g[58, g[59, g[60, g[61, g[62, g[63, g[64, g[65, g[66, g[67, g[68,
g[69]]]]]]]]]]]]]]]]]]]]]]]]]]]]]]]]]]]]]]]]]]]]]]]]]]]]]]]]]]]]]]]]]]]]]]), 16]
‘™

TagBox[
InterpretationBox[""3.127333523307754>",
3.127333523307754,
AutoDelete->True],
NumberForm[, 16] ]
NumberForm[
 4./(1 + g[0, g[1, g[2, g[3, g[4, g[5, g[6, g[7, g[8, g[9, g[10, g[11, g[12, g[13, g[14, g[15,g[16, g[17, g[18, g[19, g[20, g[21, g[22, g[23, g[24, g[25, g[26, g[27, g[28, g[29, g[30, g[31, g[32, g[33, g[34, g[35, g[36, g[37, g[38, g[39, g[40, g[41, g[42, g[43, g[44, g[45, g[46, g[47, g[48, g[49, g[50, g[51, g[52, g[53, g[54, g[55, g[56, g[57, g[58, g[59, g[60, g[61, g[62, g[63, g[64, g[65, g[66, g[67, g[68, g[69, g[70, g[71, g[72, g[73, g[74]]]]]]]]]]]]]]]]]]]]]]]]]]]]]]]]]]]]]]]]]]]]]]]]]]]]]]]]]]]]]]]]]]]]]]]]]]]), 16]
‘™

TagBox[
InterpretationBox[""3.154923023986425>",
3.1549230239864245`,
AutoDelete->True],
NumberForm[, 16 ]
NumberForm[
 4./(1 + g[0, g[2, g[3, g[4, g[5, g[6, g[7, g[8, g[9, g[10, g[11, g[12, g[13, g[14, 
g[15, g[16, g[17, g[18, g[19, g[20, g[21, g[22, g[23, g[24, g[25, g[26, g[27, g[28, g[29, g[30, g[31, g[32, g[33, g[34, g[35, g[36, g[37, g[38, g[39, g[40, g[41, g[42, g[43, g[44, g[45, g[46, g[47, g[48, g[49, g[50, g[51, g[52, g[53, g[54, g[55, g[56, g[57, g[58, g[59, g[60, g[61, g[62, g[63, g[64, g[65, g[66, g[67, g[68, g[69, g[70, g[71, g[72, g[73, g[74, g[75, g[76, g[77, g[78, g[79]]]]]]]]]]]]]]]]]]]]]]]]]]]]]]]]]]]]]]]]]]]]]]]]]]]]]]]]
]]]]]]]]]]]]]]]]]]]]]]]]), 16]
‘™

TagBox[
InterpretationBox[""3.129095094977018>",
3.1290950949770178`,
AutoDelete->True],
NumberForm[, 16] ]
NumberForm[
 4./(1 + g[0, g[1, g[2, g[3, g[4, g[5, g[6, g[7, g[8, g[9, g[10, g[11, g[12, g[13, g[14, g[15,g[16, g[17, g[18, g[19, g[20, g[21, g[22, g[23, g[24, g[25, g[26, g[27, g[28, g[29,                     g[30, g[31, g[32, g[33, g[34, g[35, g[36, g[37, g[38, g[39, g[40, g[41, g[42, g[43, g[44, g[45, g[46, g[47, g[48, g[49, g[50, g[51, g[52, g[53, g[54, g[55, g[56, g[57, g[58, g[59, g[60, g[61, g[62, g[63, g[64, g[65, g[66, g[67, g[68, g[69, g[70, g[71, g[72, g[73, g[74, g[75, g[76, g[77, g[78, g[79, g[80, g[81, g[82, g[83, g[84]]]]]]]]]]]]]]]]]]]]]]]]]]]]]]]]]]]]]]]]]]]]]]]]]]]]]]]]]]]]]]]]]]]]]]]]]]]]]]]]]]]]]), 16]
‘™

TagBox[
InterpretationBox[""3.153355324069958>",
3.153355324069958,
AutoDelete->True],
NumberForm[, 16] ]
NumberForm[
 4./(1 + g[0, g[1, g[2, g[3, g[4, g[5, g[6, g[7, g[8, g[9, g[10, g[11, g[12, g[13, g[14, g[15, g[16, g[17, g[18, g[19, g[20, g[21, g[22, g[23, g[24, g[25, g[26, g[27, g[28, g[29, g[30, g[31, g[32, g[33, g[34, g[35, g[36, g[37, g[38, g[39, g[40, g[41, g[42, g[43, g[44, g[45, g[46, g[47, g[48, g[49, g[50, g[51, g[52, g[53, g[54, g[55, g[56, g[57, g[58, g[59, g[60, g[61, g[62, g[63, g[64, g[65, g[66, g[67, g[68, g[69, g[70, g[71, g[72, g[73, g[74, g[75, g[76, g[77, g[78, g[79, g[80, g[81, g[82, g[83, g[84, g[85, g[86, g[87, g[88, g[89]]]]]]]]]]]]]]]]]]]]]]]]]]]]]]]]]]]]]]]]]]]]]]]]]]]]]]]]]]]]]]]]]]]]]]]]]]]]]]]]]]]]]]]]]]), 16]
‘™

TagBox[
InterpretationBox[""3.130483257145759>",
3.1304832571457593`,
AutoDelete->True],
NumberForm[, 16] ]
NumberForm[
 4./(1 + g[0, g[1, g[2, g[3, g[4, g[5,  g[6, g[7, g[8, g[9, g[10, g[11, g[12, g[13, g[14, g[15,                    g[16, g[17, g[18, g[19, g[20, g[21, g[22, g[23, g[24, g[25, g[26, g[27, g[28, g[29, g[30, g[31, g[32, g[33, g[34, g[35, g[36, g[37, g[38, g[39, g[40, g[41, g[42, g[43, g[44, g[45, g[46, g[47, g[48, g[49, g[50, g[51, g[52, g[53, g[54, g[55, g[56, g[57, g[58, g[59, g[60, g[61, g[62, g[63, g[64, g[65, g[66, g[67, g[68, g[69, g[70, g[71, g[72, g[73, g[74, g[75, g[76, g[77, g[78, g[79, g[80, g[81, g[82, g[83, g[84, g[85, g[86, g[87, g[88, g[89, g[90, g[91, g[92, g[93, g[94]]]]]]]]]]]]]]]]]]]]]]]]]]]]]]]]]]]]]]]]]]]]]]]]]]]]]]]]]]]]]]]]]]]]]]]]]]]]]]]]]]]]]]]]]]]]]]]), 16]
‘™

TagBox[
InterpretationBox[""3.152117511448855",
3.152117511448855,
AutoDelete->True],
NumberForm[, 16] ]
NumberForm[
 4./(1 + g[0, g[1, g[2, g[3, g[4, g[5, g[6, g[7, g[8, g[9, g[10, g[11, g[12, g[13, 
g[14, g[15, g[16, g[17, g[18, g[19, g[20, g[21, g[22, g[23, g[24, g[25, g[26, g[27, g[28, g[29, g[30, g[31, g[32, g[33, g[34, g[35, g[36, g[37, g[38, g[39, g[40, g[41, g[42, g[43, g[44, g[45, g[46, g[47, g[48, g[49, g[50, g[51, g[52, g[53, g[54, g[55, g[56, g[57, g[58, g[59, g[60, g[61, g[62, g[63, g[64, g[65, g[66, g[67, g[68, g[69, g[70, g[71, g[72, g[73, g[74, g[75, g[76, g[77, g[78, g[79, g[80, g[81, g[82, g[83, g[84, g[85, g[86, g[87, g[88, g[89, g[90, g[91, g[92, g[93, g[94, g[95, g[96, g[97, g[98, g[99]]]]]]]]]]]]]]]]]]]]]]]]]]]]]]]]]]]]]]]]]]]]]]]]]]]]]]]]]]]]]]]]]]]]]]]]]]]]]]]]]]]]]]]]]]]]]]]]]]]]), 16]
‘™

TagBox[
InterpretationBox[""3.131593903583553",
3.1315939035835534`,
AutoDelete->True],
NumberForm[, 16] ]
"Numbers 1-4 From the Function n=1-4"

NumberForm[4./(1 + g[0]), 16]
‘™

TagBox[
InterpretationBox[""2.666666666666667>",
2.6666666666666665`,
AutoDelete->True],
NumberForm[, 16] ]
NumberForm[4./(1 + g[0, g[1]]), 16]
‘™

TagBox[
InterpretationBox[""2.8>",
2.8,
AutoDelete->True],
NumberForm[, 16] ]
NumberForm[4./(1 + g[0, g[1, g[2]]]), 16]
‘™

TagBox[
InterpretationBox["3.428571428571428>",
3.4285714285714284`,
AutoDelete->True],
NumberForm[, 16] ]
NumberForm[4./(1 + g[0, g[1, g[2, g[3]]]]), 16]
‘™

TagBox[
InterpretationBox["2.911111111111111>",
2.911111111111111,
AutoDelete->True],
NumberForm[, 16] ]

"Numbers 6-9 from the Function n=6-9"

NumberForm[4./(1 + g[0, g[1, g[2, g[3, g[4, g[5]]]]]]), 16]
‘™

TagBox[
InterpretationBox[""2.980708180708181>",
2.9807081807081808`,
AutoDelete->True],
NumberForm[, 16] ]
NumberForm[4./(1 + g[0, g[1, g[2, g[3, g[4, g[5, g[6]]]]]]]), 16]
‘™

TagBox[
InterpretationBox[""3.280808080808081",
3.2808080808080806`,
AutoDelete->True],
NumberForm[, 16] ]
NumberForm[4./(1 + g[0, g[1, g[2, g[3, g[4, g[5, g[6, g[7]]]]]]]]), 16]
‘™

TagBox[
InterpretationBox[""3.280808080808081",
3.2808080808080806`,
AutoDelete->True],
NumberForm[, 16] ]
NumberForm[4./(1 + g[0, g[1, g[2, g[3, g[4, g[5, g[6, g[7]]]]]]]]), 16]
‘™

TagBox[
InterpretationBox[""3.019032601385542>",
3.0190326013855424`,
AutoDelete->True],
NumberForm[, 16]]
NumberForm[4./(1 + g[0, g[1, g[2, g[3, g[4, g[5, g[6, g[7, g[8]]]]]]]]]), 16]
‘™

\subsection{Mathematic Iteration code: Continued Fractions for $\pi$ Percentage Error}
\begin{multicols}{2}
\begin{itemize}
\item In[152]:= 2.666666666666667/π
In[153]:= 1 - 0.8488263631567753`
In[154]:= 100*0.1511736368432247`
Out[154]= 15.1174

\item InIn[155]:= 2.800000000000000/π
In[156]:= 1 - 0.8912676813146139`
In[157]:= 100*0.10873231868538613`
Out[157]= 10.8732

\item InIn[158]:= 3.428571428571428/π
In[159]:= 1 - 1.091348181201568`
In[160]:= 100*-0.09134818120156796`
Out[160]= -9.13482

\item InIn[161]:= 2.9111111111111111/π
In[162]:= 1 - 0.9266354464461461`
In[163]:= 100*0.07336455355385385`
Out[163]= 7.33646

\item InIn[164]:= 3.331601731601732/π
In[165]:= 1 - 1.0604817679958674`
In[166]:= 100*-0.06048176799586735`
Out[166]= -6.04818

\item InIn[167]:= 2.980708180708181/π
In[168]:= 1 - 0.948788881748315`
In[169]:= 100*0.051211118251685006`
Out[169]= 5.12111

\item InIn[170]:= 3.280808080808081/π
In[171]:= 1 - 1.044313646792881`
In[172]:= 100*-0.04431364679288108`
Out[172]= -4.43136

\item InIn[173]:= 3.019032601385542/π
In[174]:= 1 - 0.9609879237321854`
In[175]:= 100*0.039012076267814555`
Out[175]= 3.90121

\item InIn[176]:= 3.250989945726788/π
In[177]:= 1 - 1.0348222396089417`
In[178]:= 100*-0.03482223960894171`
Out[178]= -3.48222

\item InIn[179]:= 3.042842125195067/π
In[180]:= 1 - 0.9685667305460856`
In[181]:= 100*0.03143326945391445`
Out[181]= 3.14333

\item InIn[182]:= 3.207889026381334/π
In[183]:= 1 - 1.0211027908776735`
In[184]:= 100*-0.021102790877673527`
Out[184]= -2.11028

\item InIn[185]:= 3.091748884841698/π
In[186]:= 1 - 0.9841342356428227`
In[187]:= 100*0.015865764357177348`
Out[187]= 1.58658

\item InIn[188]:= 3.181512659824787/π
In[189]:= 1 - 1.012706932641117`
In[190]:= 100*-0.012706932641117108`
Out[190]= -1.27069

\item InIn[191]:= 3.108305614026907/π
In[192]:= 1 - 0.9894044062253424`
In[193]:= 100*0.010595593774657575`
Out[193]= 1.05956

\item InIn[194]:= 3.170134928816534/π
In[195]:= 1 - 1.0090852883788504`
In[196]:= 100*-0.009085288378850365`
Out[196]= -0.908529

\item InIn[197]:= 3.116612184235621/π
In[198]:= 1 - 0.9920484696430558`
In[199]:= 100*0.00795153035694418`
Out[199]= 0.795153

\item InIn[200]:= 3.163801158726882/π
In[201]:= 1 - 1.007069186742499`
In[202]:= 100*-0.007069186742499012`
Out[202]= -0.706919

\item InIn[203]:= 3.121602653391091/π
In[204]:= 1 - 0.9936369853119372`
In[205]:= 100*0.006363014688062774`
Out[205]= 0.636301

\item InIn[206]:= 3.123742628422234/π
In[207]:= 1 - 0.9943181605205365`
In[208]:= 100*0.005681839479463546`
Out[208]= 0.568184

\item InIn[209]:= 3.124931773880152/π
In[210]:= 1 - 0.9946966772759024`
In[211]:= 100*0.005303322724097614`
Out[211]= 0.530332

\item InIn[212]:= 3.156972717366711/π
In[213]:= 1 - 1.0048956263503301`
In[214]:= 100*-0.00489562635033014`
Out[214]= -0.489563

\item InIn[215]:= 3.127333523307754/π
In[216]:= 1 - 0.9954611778628444`
In[217]:= 100*0.004538822137155618`
Out[217]= 0.453882

\item InIn[218]:= 3.154923023986425/π
In[219]:= 1 - 1.0042431886837397`
In[220]:= 100*-0.00424318868373974`
Out[220]= -0.424319

\item InIn[221]:= 3.129095094977018/π
In[222]:= 1 - 0.9960219035403924`
In[223]:= 100*0.003978096459607561`
Out[223]= 0.39781

\item InIn[224]:= 3.153355324069958/π
In[225]:= 1 - 1.0037441743017588`
In[226]:= 100*-0.0037441743017587736`
Out[226]= -0.374417

\item InIn[227]:= 3.130483257145759/π
In[228]:= 1 - 0.9964637692823288`
In[229]:= 100*0.003536230717671174`
Out[229]= 0.353623

\item InIn[230]:= 3.152117511448855/π
In[231]:= 1 - 1.0033501663072186`
In[232]:= 100*-0.003350166307218627`
Out[232]= -0.335017

\item InIn[233]:= 3.131593903583553/π
In[234]:= 1 - 0.9968172990235336`
In[235]:= 100*0.003182700976466446`
Out[235]= 0.31827
\end{itemize}
\end{multicols}
\newpage

\section{Mathematic Iteration code: Viete's Nested Radicals Formula for $\pi$}
"Viete Formula"
$f[x_] := Sqrt[2 + x]; p[n_] := 2^(n + 1) Sqrt[2 - Nest[f, Sqrt[2], n - 1]];$ 
$a[n_] := IntegerPart[3.3 n];$
For[i = 0, i <= 100, i = i + 5, 
 Print[i " SymbA " Text[Style[Row[[N[p[i], a[i]]]]]]]]For[i = 0, i <= 10, i = i + 1,  Print[i " SymbA " Text[Style[Row[{N[p[i], a[i]]}]]]]]

1 " SymbA " 3.061467458921242

2 " SymbA " 3.121445152263491

3 " SymbA " 3.136548490541725

4 " SymbA " 3.140331156952385

5 " SymbA " 3.141277250932773

6 " SymbA " 3.141513801175428

7 " SymbA " 3.141572940255612

8 " SymbA " 3.141587725373528

9 " SymbA " 3.141591413562714

10 " SymbA " 3.141592345570118

15 " SymbA " 3.141592653288993

20 " SymbA " 3.141592653589499

25 " SymbA " 3.141592653589793

30 " SymbA " 3.141592653589793

35 " SymbA " 3.141592653589793

40 " SymbA " 3.141592653589793

45 " SymbA " 3.141592653589793

50 " SymbA " 3.141592653589793

55 " SymbA " 3.141592653589793

60 " SymbA " 3.141592653589793

65 " SymbA " 3.141592653589793

70 " SymbA " 3.141592653589793

75 " SymbA " 3.141592653589793

80 " SymbA " 3.141592653589793

85 " SymbA " 3.141592653589793

90 " SymbA " 3.141592653589793

95 " SymbA " 3.141592653589793

100 " SymbA" 3.141592653589793
\subsection{Mathematic Iteration code: Viete's Nested Radicals Percentage Error}
\begin{multicols}{2}
\begin{itemize}
\item 3.061467458921242/π
1 - 0.9744953584045994`
100*0.025504641595400557`
2.55046

\item 3.121445152263491/π
1 - 0.993586851145937`
100*0.006413148854062967`
0.641315

\item 3.136548490541725/π
1 - 0.9983943930342769`
100*0.0016056069657230942`
0.160561

\item 3.140331156952385/π
1 - 0.9995984531489255`
100*0.00040154685107451904`
0.0401547

\item 3.141277250932773/π
1 - 0.9998996042161419`
100*0.00010039578385812042`
0.0100396

\item 3.141513801175428/π
1 - 0.9999749004969581`
100*0.000025099503041858817`
0.00250995

\item 3.141572940255612/π
1 - 0.9999937250508405`
100*6.274949159501553`*-6
0.000627495

\item 3.141587725373528/π
1 - 0.9999984313000415`
100*1.5686999584874073`*-6
0.00015687

\item 3.141591413562714/π
1 - 0.9999996052871216`
100*3.9471287838210856`*-7
0.0000394713
\end{itemize}
\end{multicols}

\newpage
\section{Mathematic Iteration code: Zeta (2),(4),(6),(8)}

N[Sqrt[6*Sum[1/k2, {k, 1, 5}]], 15]
3.09466952411370
N[Sqrt[6*Sum[1/k2, {k, 1, 10}]], 15]
3.04936163598207
N[Sqrt[6*Sum[1/k2, {k, 1, 15}]], 15]
3.07938982603209
N[Sqrt[6*Sum[1/k2, {k, 1, 20}]], 15]
3.09466952411370
N[Sqrt[6*Sum[1/k2, {k, 1, 25}]], 15]
3.10392339170058
N[Sqrt[6*Sum[1/k2, {k, 1, 30}]], 15]
3.11012872814126
N[Sqrt[6*Sum[1/k2, {k, 1, 35}]], 15]
3.11457886229313
N[Sqrt[6*Sum[1/k2, {k, 1, 40}]], 15]
3.11792619829938
N[Sqrt[6*Sum[1/k2, {k, 1, 45}]], 15]
3.12053546308708
N[Sqrt[6*Sum[1/k2, {k, 1, 50}]], 15]
3.12262652293373
N[Sqrt[6*Sum[1/k2, {k, 1, 55}]], 15]
3.12433980504914
N[Sqrt[6*Sum[1/k2, {k, 1, 60}]], 15]
3.12576920214052
N[Sqrt[6*Sum[1/k2, {k, 1, 65}]], 15]
3.12697987310384
N[Sqrt[6*Sum[1/k2, {k, 1, 70}]], 15]
3.12801845342065
N[Sqrt[6*Sum[1/k2, {k, 1, 75}]], 15]
3.12891920064047
N[Sqrt[6*Sum[1/k2, {k, 1, 80}]], 15]
3.12970784547462
N[Sqrt[6*Sum[1/k2, {k, 1, 85}]], 15]
3.13040408931831
N[Sqrt[6*Sum[1/k2, {k, 1, 90}]], 15]
3.13102327252367
N[Sqrt[6*Sum[1/k2, {k, 1, 95}]], 15]
3.13157751780151
N[Sqrt[6*Sum[1/k2, {k, 1, 100}]], 15]
3.13207653180911 \\ \\

N[Surd[90*Sum[1/k4, {k, 1, 5}], 4], 15]
3.14016117947426
N[Surd[90*Sum[1/k4, {k, 1, 10}], 4], 15]
3.14138462246697
N[Surd[90*Sum[1/k4, {k, 1, 15}], 4], 15]
3.14152783068467
N[Surd[90*Sum[1/k4, {k, 1, 20}], 4], 15]
3.14156460959141
N[Surd[90*Sum[1/k4, {k, 1, 25}], 4], 15]
3.14157807684660
N[Surd[90*Sum[1/k4, {k, 1, 30}], 4], 15]
3.14158413278489
N[Surd[90*Sum[1/k4, {k, 1, 35}], 4], 15]
3.14158724909022
N[Surd[90*Sum[1/k4, {k, 1, 40}], 4], 15]
3.14158901347572
N[Surd[90*Sum[1/k4, {k, 1, 45}], 4], 15]
3.14159008631043
N[Surd[90*Sum[1/k4, {k, 1, 50}], 4], 15]
3.14159077577492
N[Surd[90*Sum[1/k4, {k, 1, 55}], 4], 15]
3.14159123889581
N[Surd[90*Sum[1/k4, {k, 1, 60}], 4], 15]
3.14159156142938
N[Surd[90*Sum[1/k4, {k, 1, 65}], 4], 15]
3.14159179291850
N[Surd[90*Sum[1/k4, {k, 1, 70}], 4], 15]
3.14159196334879
N[Surd[90*Sum[1/k4, {k, 1, 75}], 4], 15]
3.14159209159432
N[Surd[90*Sum[1/k4, {k, 1, 80}], 4], 15]
3.14159218993944
N[Surd[90*Sum[1/k4, {k, 1, 85}], 4], 15]
3.14159226661411
N[Surd[90*Sum[1/k4, {k, 1, 90}], 4], 15]
3.14159232727297
N[Surd[90*Sum[1/k4, {k, 1, 95}], 4], 15]
3.14159237588858
N[Surd[90*Sum[1/k4, {k, 1, 100}], 4], 15]
3.14159241530737 \\ \\

N[Surd[945*Sum[1/k6, {k, 1, 5}], 6], 15]
3.14157300346359
N[Surd[945*Sum[1/k6, {k, 1, 10}], 6], 15]
3.14159185608168
N[Surd[945*Sum[1/k6, {k, 1, 15}], 6], 15]
3.14159253913011
N[Surd[945*Sum[1/k6, {k, 1, 20}], 6], 15]
3.14159262524305
N[Surd[945*Sum[1/k6, {k, 1, 25}], 6], 15]
3.14159264406125
N[Surd[945*Sum[1/k6, {k, 1, 30}], 6], 15]
3.14159264969505
N[Surd[945*Sum[1/k6, {k, 1, 35}], 6], 15]
3.14159265176594
N[Surd[945*Sum[1/k6, {k, 1, 40}], 6], 15]
3.14159265264583
N[Surd[945*Sum[1/k6, {k, 1, 45}], 6], 15]
3.14159265306227
N[Surd[945*Sum[1/k6, {k, 1, 50}], 6], 15]
3.14159265327654
N[Surd[945*Sum[1/k6, {k, 1, 55}], 6], 15
3.14159265339440
N[Surd[945*Sum[1/k6, {k, 1, 60}], 6], 15]
3.14159265346284
N[Surd[945*Sum[1/k6, {k, 1, 65}], 6], 15]
3.14159265350444
N[Surd[945*Sum[1/k6, {k, 1, 70}], 6], 15]
3.14159265353070
N[Surd[945*Sum[1/k6, {k, 1, 75}], 6], 15]
3.14159265354784
N[Surd[945*Sum[1/k6, {k, 1, 80}], 6], 15]
3.14159265355935
N[Surd[945*Sum[1/k6, {k, 1, 85}], 6], 15]
3.14159265356727
N[Surd[945*Sum[1/k6, {k, 1, 90}], 6], 15]
3.14159265357284
N[Surd[945*Sum[1/k6, {k, 1, 95}], 6], 15]
3.14159265357684
N[Surd[945*Sum[1/k6, {k, 1, 100}], 6], 15]
3.14159265357975 \\ \\

N[Surd[9450*Sum[1/k8, {k, 1, 5}], 8], 15]
3.14159231269578
N[Surd[9450*Sum[1/k8, {k, 1, 10}], 8], 15]
3.14159264970117
N[Surd[9450*Sum[1/k8, {k, 1, 15}], 8], 15]
3.14159265333235
N[Surd[9450*Sum[1/k8, {k, 1, 20}], 8], 15]
3.14159265355327
N[Surd[9450*Sum[1/k8, {k, 1, 25}], 8], 15]
3.14159265358185
N[Surd[9450*Sum[1/k8, {k, 1, 30}], 8], 15]
3.14159265358752
N[Surd[9450*Sum[1/k8, {k, 1, 35}], 8], 15]
3.14159265358901
N[Surd[9450*Sum[1/k8, {k, 1, 40}], 8], 15]
3.14159265358948
N[Surd[9450*Sum[1/k8, {k, 1, 45}], 8], 15]
3.14159265358966
N[Surd[9450*Sum[1/k8, {k, 1, 50}], 8], 15]
3.14159265358973
N[Surd[9450*Sum[1/k8, {k, 1, 55}], 8], 15]
3.14159265358976
N[Surd[9450*Sum[1/k8, {k, 1, 60}], 8], 15]
3.14159265358977
N[Surd[9450*Sum[1/k8, {k, 1, 65}], 8], 15]
3.14159265358978
N[Surd[9450*Sum[1/k8, {k, 1, 70}], 8], 15]
3.14159265358979
N[Surd[9450*Sum[1/k8, {k, 1, 75}], 8], 15]
3.14159265358979
N[Surd[9450*Sum[1/k8, {k, 1, 80}], 8], 15]
3.14159265358979
N[Surd[9450*Sum[1/k8, {k, 1, 85}], 8], 15]
3.14159265358979
N[Surd[9450*Sum[1/k8, {k, 1, 90}], 8], 15]
3.14159265358979
N[Surd[9450*Sum[1/k8, {k, 1, 95}], 8], 15]
3.14159265358979
N[Surd[9450*Sum[1/k8, {k, 1, 100}], 8], 15]
3.14159265358979

\subsection{Mathematic Iteration code:  Zeta (2),(4),(6),(8)}
\begin{multicols}{2}
\begin{itemize}
\item 3.09466952411372/ 3.141592653589793
1 - 0.985063903997084
0.01493609600291601`*100
1.49361

\item 3.04936163598207/3.141592653589793
1 - 0.9706419552826705`
0.0293580447173295`*100
2.9358

\item 3.07938982603209/3.141592653589793
1 - 0.9802002250397976`
100*0.019799774960202354`
1.97998

\item 3.09466952411378/3.141592653589793
1 - 0.985063903997103`
100*0.014936096002897026`
1.49361

\item 3.10392339170058/3.141592653589793
1 - 0.9880095015354172`
100*0.011990498464582777`
1.19905

\item 3.11012872814126/3.141592653589793
1 - 0.9899847214715822`
100*0.01001527852841777`
1.00153

\item 3.11457886229313/3.141592653589793
1 - 0.9914012431669665`
100*0.00859875683303346`
0.859876

\item 3.11792619829938/3.141592653589793
1 - 0.9924667333101348`
100*0.007533266689865203`
0.753327

\item 3.12053546308708/3.141592653589793
1 - 0.993297288087731`
100*0.00670271191226901`
0.670271

\item 3.12262652293373/3.141592653589793
1 - 0.9939628931095217`
100*0.0060371068904783165`
0.603711

\item 3.12433980504914/3.141592653589793
1 - 0.9945082477446786`
100*0.0054917522553213916`
0.549175

\item 3.12576920214052/3.141592653589793
1 - 0.9949632389701472`
100*0.005036761029852843`
0.503676

\item 3.12697987310384/3.141592653589793
1 - 0.9953486075066875`
100*0.004651392493312478`
0.465139

\item 3.12801845342065/3.141592653589793
1 - 0.9956791978891241`
100*0.004320802110875932`
0.43208

\item 3.12891920064047/3.141592653589793
1 - 0.9959659146341453`
100*0.0040340853658547005`
0.403409

\item 3.12970784547462/3.141592653589793
1 - 0.996216948081543`
100*0.003783051918456959`
0.378305

\item 3.13040408931831/3.141592653589793
1 - 0.9964385693801843`
100*0.003561430619815731`
0.356143

\item 3.13102327252367/3.141592653589793
1 - 0.9966356615158093`
100*0.00336433848419071`
0.336434

\item 3.13157751780151/3.141592653589793
1 - 0.9968120832671165`
\item 100*0.0031879167328835445`
0.318792

\item 3.13207653180911/3.141592653589793
1 - 0.9969709243590797`
100*0.0030290756409202535`
0.302908
\end{itemize}
\end{multicols}


\newpage
\section{Bibliography}
\begin{thebibliography}{12}
\bibitem{cite_key1} BECKMANN, Petr. $1971$. Viete’s formula in the history of pi [online]. [Accessed $19$ August $2013$]. Wikipedia. Available from World Wide Web: http://en.wikipedia.org/wiki/Viete'sformula
\\
\bibitem{cite_key2} FIRK, Frank. 2001. A proof of the Gregory-Leibniz series and new series for calculating π. Arxiv. 1(1), pp. 3-4
\\
\bibitem{cite_key3} GOUREVITCH, Boris. 2005. Newton’s Formula [online]. [Acessed 19th August 2013]. Pi314. Available from World Wide Web: http://www.pi314.net/eng/newton.php
\\
\bibitem{cite_key4} GUSMORINO, P. 1998. Gregory, Leibniz, and Machin. Using the Arctangent Infinite Series. 12(1). pp.3-4.
\\
\bibitem{cite_key5} HAYASHI, Ko. 1996. Infinite Expressions for Pi [online]. [Accessed 19 August 2013]. Geom. Available from World Wide Web: http://www.geom.uiuc.edu/~huberty/math5337/groupe/expresspi.html
\\
\bibitem{cite_key6} Hazewinkel, Michiel. 2001. Parseval Equality [online]. [Accessed 19 August 2013]. Wikipedia. Available from World Wide Web: http://en.wikipedia.org/wiki/Parseval'sidentity
\\
\bibitem{cite_key7} JONATHAN, Borwein. 2013. Chronology of computation of π [online]. [Accessed 19 August 2013]. Wikipedia. Available from World Wide Web:
http://en.wikipedia.org/wiki/Chronologyofcomputationofπ
\\
\bibitem{cite_key8} JONES, William. 1980. Continued Fractions [online]. [Accessed 19 August 2013]. Wikipedia. Available from World Wide Web: http://en.wikipedia.org/wiki/Continuedfraction
\\
\bibitem{cite_key9} Lamb, Evelyn. 2013. How Much Pi Do You Need? [online]. [Accessed 19 August 2013]. ScientificeAmerican. Available from World Wide Web: http://blogs.scientificamerican.com/observations/2012/07/21/how-much-pi-do-you-need/
\\
\bibitem{cite_key10} LEJEUNE,Dirichlet. 1829. Fourier’s formula for 2pi-periodic functions using sines and cosines. Wikipedia. Available from World Wide Web: http://en.wikipedia.org/wiki/Fourierseries
\\
\bibitem{cite_key11} LORENTZEN, Lisa. 1992. Continued fractions with Applications [online]. [Accessed 19 August 2013]. Wikipedia. Available from World Wide Web:
http://en.wikipedia.org/wiki/Generalizedcontinuedfraction
\\
\bibitem{cite_key12} LYNN, Ben. 2006. The Wallis Product [online]. [Accessed 19 August 2013]. Stanford. Crypto. Available from World Wide Web: http://crypto.stanford.edu/pbc//notes/pi/wallis.html
\end{thebibliography}
\end{document}